\documentclass[a4paper,12 pt,leqno]{article}
\usepackage{amsmath,amsfonts,here}
\font\tencmmib=cmmib10
\skewchar\tencmmib='177
\font\sevencmmib=cmmib7
\skewchar\sevencmmib='177
\font\fivecmmib=cmmib5
\skewchar\fivecmmib='177
\newfam\cmmibfam
\textfont\cmmibfam=\tencmmib
\scriptfont\cmmibfam=\sevencmmib
\scriptscriptfont\cmmibfam=\fivecmmib
\font\tencmbsy=cmbsy10
\skewchar\tencmbsy='60
\font\sevencmbsy=cmbsy7
\skewchar\sevencmbsy='60
\font\fivecmbsy=cmbsy5
\skewchar\fivecmbsy='60
\newfam\cmbsyfam
\textfont\cmbsyfam=\tencmbsy
\scriptfont\cmbsyfam=\sevencmbsy
\scriptscriptfont\cmbsyfam=\fivecmbsy

\begin{document}
\newtheorem{theo}{Theorem}[section]
\newtheorem{lemme}[theo]{Lemma}
\newtheorem{cor}[theo]{Corollary}
\newtheorem{defi}[theo]{Definition}
\newtheorem{prop}[theo]{Proposition}
\newtheorem{problem}[theo]{Problem}
\newcommand{\beq}{\begin{eqnarray}}
\newcommand{\enq}{\end{eqnarray}}
\newcommand{\be}{\begin{eqnarray*}}
\newcommand{\en}{\end{eqnarray*}}
\newcommand{\Td}{\mathbb T^d}
\newcommand{\T}{\mathbb T}
\newcommand{\Rd}{\mathbb R^d}
\newcommand{\R}{\mathbb R}
\newcommand{\Zd}{\mathbb Z^d}
\newcommand{\Z}{\mathbb Z}
\newcommand{\Linf}{L^{\infty}}
\newcommand{\dt}{\partial_t}
\newcommand{\Dt}{\frac{d}{dt}}
\newcommand{\Dtt}{\frac{d^2}{dt^2}}
\newcommand{\demi}{\frac{1}{2}}
\newcommand{\vf}{\varphi}
\newcommand{\epu}{_{\epsilon}}
\newcommand{\ep}{^{\epsilon}}
\newcommand{\bfi}{{\mathbf \Phi}}
\newcommand{\bpsi}{{\mathbf \Psi}}
\newcommand{\ds}{\displaystyle}
\newcommand{\bx}{{\mathbf X}}
\newcommand{\bg}{{\mathbf g}}
\newcommand{\bv}{{\mathbf v}}
\title{
On the regularity of the polar factorization for time dependent maps
}
\maketitle 
\bibliographystyle{plain}

\begin{center}
G. Loeper\footnotemark[1]\footnotemark[2]
\end{center}
\footnotetext[1]{Universit\'e de
Nice-Sophia-Antipolis 
\hfill  email:loeper@math.unice.fr}
\footnotetext[2]{Fields Institute \& University of Toronto}
\begin{abstract}
We consider the polar factorization of vector valued mappings, introduced in \cite{Br1}, in the 
case of a family of mappings depending on a parameter. We investigate
the regularity with respect to this parameter of the terms of the polar factorization  
by constructing some a priori bounds.
To do so, we consider the linearization of the associated Monge-Amp\`ere equation.
\end{abstract}  
\section{Introduction}


\subsubsection*{Polar factorization and Monge-Amp\`ere equation}
Brenier in \cite{Br1} showed that given $\Omega$ a bounded open set of $\Rd$
such that $|\partial\Omega|=0$, with $|.|$ the Lebesgue measure of $\Rd$,
every Lebesgue measurable 
mapping ${\bf X}\in L^2(\Omega, \Rd)$  satisfying the 
non-degeneracy condition 
\beq\label{2N}
\forall B\subset\Rd \textrm{ measurable, } |B|=0 \Rightarrow  |{\bf X}^{-1}(B)|=0
\enq
can be factorized in the following (unique) way:
\beq\label{2polar}
{\bf X}=\nabla\bfi\circ {\bf g},
\enq
where $\bfi$ is a convex function and ${\bf g}$ belongs to $G(\Omega)$ the set of Lebesgue-measure preserving mappings of $\Omega$, defined by
\beq\label{2meas-pres}
{\bf g}  \in G(\Omega) \iff
\forall f \in C_b(\Omega), \int_{\Omega}f({\bf g}(x)) \ dx =\int_{\Omega}f(x) \ dx, 
\enq
where $C_b$ is the set of bounded continuous functions.
If $da$ denotes the Lebesgue measure of $\Omega$, the push-forward of  $da$ by ${\bf X}$, that we denote $\bx\# da$, is the measure $\rho$ defined by
\beq\label{2defrho}
\forall f\in C_b(\Rd),\; \int_{\Rd} f d\rho=\int_{\Omega} f({\bf X}(a))da. 
\enq
One sees first that the condition (\ref{2N}) is equivalent to the fact that $\rho$ is absolutely continuous with respect to the Lebesgue measure, or has a density in $L^1(\Rd, dx)$. 
Then $\bfi$ satisfies in $\Omega$ the
Monge-Amp\`ere equation: 
\be 
\rho(\nabla\bfi(x))\det D^2\bfi(x)=1 
\en
in the following weak sense:
\beq\label{2defPhi}
\forall g\in C_b(\Rd), \int_{\Omega}g(\nabla\bfi(y))dy =\int_{\Rd} g(x)d\rho(x).
\enq
$\bpsi$, the Legendre transform of $\bfi$, defined by
\beq\label{2deflegendre}
\bpsi(y)=\sup_{x\in \Omega} \{x\cdot y -\bfi(x)\},
\enq
satisfies  the Monge-Amp\`ere equation 
\be
\det D^2\bpsi(x)=\rho(x) 
\en 
in the following weak sense:
\beq
&&\forall f\in C_b(\Omega), \int_{\Rd} f(\nabla\bpsi(x))d\rho(x)=\int_\Omega f(y)dy\label{2defPsi}.
\enq
Note that the existence and uniqueness of the pair $\nabla\bfi,\nabla\bpsi$ 
and the validity of (\ref{2defPhi}) is not subject to the condition (\ref{2N}) (see \cite{Vi} Th 2.12 for this precise fact, and for a complete reference on polar factorization and optimal transportation). 
However (\ref{2defPsi}) may not hold.
Note also that this formulation of the second boundary value problem
for the Monge-Amp\`ere equation is strictly weaker than the Aleksandrov formulation (see \cite{Ca3} where the different formulations are compared and where
it is shown that they may not co\"incide if some extra conditions are not satisfied).

\vspace{1cm}
{\bf The periodic  case}
The polar factorization of maps on general Riemannian manifolds has been treated by \cite{Mc2}, and also in the  particular case of the flat torus by \cite{Co}.  Given $\bx$ a mapping of $\Td=\Rd/\Zd$ into itself, 
we look for a pair $(\bfi,\bg)$ such that
\begin{enumerate}
\item $\bg$ is measure preserving from $\Td$ into itself,
\item $\bfi$ is convex from $\Rd$ to $\R$ and $\bfi -|x|^2/2$ is periodic,
\item $\bx = \nabla\bfi \circ \bg$ (Note that the condition above ensures that
$\nabla\bfi-x$ is $\Zd$ periodic).
\end{enumerate}
Then under the non-degeneracy condition (\ref{2N}), there exists a unique such pair $(\bg, \nabla\bfi)$.

\vspace{1cm}

\subsubsection*{Introducing the time-dependence}
In this paper we are interested in the following problem: given a ``time'' dependent family of mappings $t\to {\bf X}(t,.)$, where for all $t$, $\bx(t)$ maps $\Omega$ in $\Rd$, 
we investigate the regularity  of the curve $t\rightarrow ({\bf g}(t,.), \bfi(t,.), \bpsi(t,.))$.

 We state different results under different assumptions.
The weakest assumption is that $\rho={\bf X}{\#} da$, $\bx$ and $\dt {\bf X}$  belong to $\Linf$ in time and space.
In this case $\dt \nabla \bfi$ and $\dt \bg$ are bounded as measures (Th. \ref{2main1}). 

Under the additional assumption that $\rho$ is close to $1$ (or actually to a continuous positive function) in $\Linf$ norm (but we do not ask for continuity), we obtain that
$\dt \bfi$ belongs to $C^{\alpha}$ for some $\alpha>0$ (Th. \ref{2main2}).
To this purpose we use a local maximum principle for solutions 
of degenerate elliptic equations (Theorem \ref{2main5}, Theorem \ref{2main4}) obtained by Murthy and Stampacchia (\cite{MuSt}) and Trudinger (\cite{Tr}), and use a result  by Caffarelli and Gutierrez (\cite{CaGu}) that establishes the  Harnack inequality for solutions of the homogeneous linearized Monge-Amp\`ere equation (Theorem \ref{2theocagu}).

The polar factorization has the following  geometrical interpretation: 
if ${\bf X} =\nabla\bfi\circ {\bf g}$, as in (\ref{2polar}), then ${\bf g}$ is the projection, in the $L^2(\Omega,\Rd)$ sense, of ${\bf X}$ on $G(\Omega)$, the set of Lebesgue measure preserving mappings. Therefore our study amounts to examine the continuity and the differentiability of the projection operator on $G(\Omega)$. We also briefly discuss
a variant of the Hodge decomposition of vector fields that appears naturally in 
this study.

Our results have an immediate application to 
the semi-geostrophic equations, a system  
arising in meteorology to model frontogenesis (see \cite{CP}). They allow in particular
to define the velocity in the physical space, a fact that was not known for weak solutions.
We discuss this application in a more extensive way in section \ref{2SG}. 

\vspace{1cm}

\subsection{Heuristics}
We present here some formal computations, assuming that all the terms considered are smooth enough. 
Suppose that $\Omega$ is bounded,
and for any $t$ we  denote by $d\rho(t,\cdot)={\bf X}(t,\cdot)\#da$ (with $da$ the
Lebesgue measure on $\Omega$) the measure defined by (\ref{2defrho}).
Then for all $t$, $\bfi(t,\cdot),\bpsi(t,\cdot)$ are as in (\ref{2defPhi},\ref{2defPsi}). 
\subsubsection*{Parallel with the Hodge decomposition of vector fields}
By differentiating (\ref{2polar}) with respect to time one finds 
\be
\dt {\bf X}(t,a)= \dt\nabla\bfi(t,{\bf g}(t,a))+D^2\bfi(t,{\bf g}(t,a))\dt {\bf g}(t,a).
\en
If ${\bf X}$ is invertible,  one can 
write 
\beq\label{2vexiste}
\dt {\bf X}(t,a)=v(t,{\bf X}(t,a))
\enq
for some ``Eulerian'' vector field $v(t,x)$ defined $d\rho$ a.e.
Note that $\rho=\bx\#da$ and $v$ will be linked through the mass conservation constraint
\beq\label{2hodge-conti}
\dt\rho + \nabla\cdot(\rho v)=0.
\enq
${\bf g}$ will then also be invertible and composing with ${\bf g}^{-1}$ one gets:
\beq\label{2hodge-ampere}
v(t,\nabla\bfi(t,x))=\dt \nabla \bfi (t,x)+D^2\bfi(t,x) w(t,x)
\enq
with $w=\dt \bg(t,\bg^{-1}(t,x))$. Since for all $t$, $\bg(t)\in G(\Omega)$, 
it follows that $w$ is divergence free.
Composing with $\nabla\bpsi=\nabla\bfi^{-1}$ we obtain
\be
v=\dt\nabla\bfi(\nabla\bpsi) + D^2\bfi  \ w(\nabla\bpsi).
\en
It is easily checked that $\bar w = D^2\bfi  \ w(\nabla\bpsi)$ satisfies
\be
\nabla\cdot(\rho \bar w)=0,
\en
therefore the second term in the decomposition (\ref{2hodge-ampere}) does not move mass. It plays the role of a divergence free vector field for a uniform density. 
\\
Note that a similar decomposition is performed in the study of the incompressible inhomogeneous Navier-Stokes equation in \cite{PLL} where 
for a given velocity field $v$, and a density $\rho>0$, one seeks to decompose $v$ as 
\be
v= \frac{1}{\rho}\nabla p + w, \hspace{.5cm} \nabla\cdot w=0. 
\en
The next proposition shows that, in the non-degenerate case where $\bfi$ is smooth and strictly convex,  the decomposition (\ref{2hodge-ampere}) is defined in an unique way.


\begin{prop}\label{2hodge}
Let $v\in L^2(\Rd, d\rho; \Rd)$, let $\bfi:\bar\Omega \to \Rd$ be $C^2$ and strictly convex on $\bar\Omega$, with $\rho= \nabla\bfi{\#} da$.
Then there exists a  unique decomposition of $v$ such that
\beq\label{2propodecompo}
v(\nabla\bfi)=\nabla p + D^2\bfi\, w
\enq
with 
$(\nabla p, w) \in L^2(\Omega;\Rd)$, 
$\nabla\cdot w=0$, $w\cdot \partial\Omega =0$.
\\
\end{prop}
{\it Proof:} We only sketch the proof of this classical result. $w$ can be found by looking for
\beq\label{2inf}
\inf_{\scriptsize\begin{array}{lll}w\in L^2(\Omega;\Rd)\\ w\cdot \partial\Omega=0\\ \nabla\cdot w =0\end{array}}\{ \int \demi w^t\cdot D^2\bfi\cdot w - v(\nabla\bfi)\cdot w\}.
\enq
Using the strict convexity of $\bfi$ we have $D^2\bfi \geq \lambda I$ on $\bar\Omega$, and we obtain that
\be
\|w\|_{L^2(\Omega)} \leq \frac{2}{\lambda}\left[\int \rho |v|^2\right]^{1/2}.
\en
The functional to minimize is strictly convex, and weakly lower semi continuous, therefore
the problem admits a minimizer. 
For the uniqueness of the decomposition, notice that if
\be
0=\nabla p +  D^2\bfi w
\en
for $\nabla p, w \in L^2$, multiplying by $w$ and integrating over $\Omega$, we get that $\nabla p, w=0$.
Therefore, if $v$ governs the evolution of $\rho$ through the equation (\ref{2hodge-conti}),
the decomposition (\ref{2propodecompo}) will co\"incide with (\ref{2hodge-ampere}) and will yield $\nabla p = \dt \nabla\bfi$.

\vspace{1cm}

\subsubsection*{The associated elliptic problems: The linearized Monge-Amp\`ere equation}
Multiplying (\ref{2hodge-ampere}) by $D^2\bfi^{-1}$, we find that $\dt\bfi$ will be solution of the following elliptic problem:
\be
\nabla\cdot (D^2\bfi^{-1}\nabla \dt\bfi) =\nabla\cdot (D^2\bfi^{-1} v(\nabla\bfi)).
\en
On the other hand, 
$\bpsi=\bfi^*$ (see (\ref{2deflegendre})) solves formally the equation 
\be
\det D^2\bpsi =\rho.
\en
Then for any $(d \times d)$ matrices $A,B$ we have
\be
\det(A+tB)=\det A + t\  \text{trace}(^tA^* B) + o(t)
\en where $A^*$ is the matrix of cofactors (or co-matrix) of $A$
and thus, formally, $\dt \bpsi$
solves the elliptic equation 
\beq
&&M_{ij}\partial_{ij}\dt \bpsi = \dt \rho\label{2malin0},
\nonumber
\enq
where $(M_{ij})_{i,j\in [1..d]}$ is the co-matrix of $D^2\bpsi$, 
given by 
\be
M=\det D^2\bpsi [D^2\bpsi]^{-1}=\rho D^2\bfi (\nabla\bpsi).
\en
Then if $M$ is the co-matrix of a second derivative matrix, for all $j\in [1..d]$
\be
\sum_{i=1}^d \partial_i M_{ij}(x) \equiv 0,
\en
and using this and the equation (\ref{2hodge-conti}), we obtain  a divergence formulation of the problem:
\beq\label{2diverg}
\nabla\cdot(M \nabla\partial_t\bpsi)=\dt \rho=-\nabla\cdot(\rho v).
\enq
In the case where $\rho$ is smooth and supported in a convex set, it will be shown using classical elliptic regularity  and results on Monge-Amp\`ere equation, that the decomposition holds (Proposition \ref{2main3a}) and that the terms are smooth. 
\\
For a generic, non-necessarily smooth $\rho$, 
we see that the difficulty will be coming from the lack of regularity and 
ellipticity of this equation. Indeed we only know a-priori that $D^2\bfi$ is a measure.
If $\rho$ is close to 1 in $\Linf$ norm, we get that $D^2\bfi$ is in $L^p_{loc}$ for some $p<\infty$,  and thus non necessarily uniformly elliptic.

\section{Results}
\subsubsection*{Notations}
In the remainder of the paper $\Omega$ will be kept fixed once for all and chosen bounded and convex. We will furthermore assume for simplicity (although one may possibly remove this assumption  through approximation) that
it is  smooth and strictly convex. 
\\
The Lebesgue measure of $\Omega$, $\chi_{\Omega}{\cal L}^d$, will be denoted in short $da$.
\\
For compatibility $\rho$ will be a probability measure on $\Rd$ and $\Omega$ of Lebesgue measure one.
\\
${\cal M}(\Omega)$ will design the set of (possibly vector valued) bounded measures on $\Omega$, with norm
$\|.\|_{{\cal M}(\Omega)}$.
\\
For $M$ a $(d\times d)$ matrix, and $u,v$ two vectors of $\Rd$, $u M v$ will denote $\sum_{i,j} u_iM_{ij} v_j$.
\\
$I$ will be an non-empty open interval of $\R$.
\\
We still use $d\rho(t,\cdot)=\bx(t,\cdot)\#da$, 
the functions $\bfi(t,\cdot), \bpsi(t,\cdot)$ will be as in (\ref{2defPhi}, \ref{2defPsi}) with $(\rho(t,\cdot), \Omega)$.
Since they are defined only up to a constant, we will impose the 
condition:
\beq\label{2jauge}
\forall t\in I, \ \int_{\Omega}\bfi(t,x) \ dx=0,
\enq
and this sets also $\bpsi$ through the relation $\bpsi=\bfi^*$.

\vspace{1cm}
\begin{theo}\label{2main1}
Let $\Omega,I$ be as above, let ${\bf X}: I \times \Omega \rightarrow\Rd.$ Let, for any $t\in I$, $d\rho(t,\cdot)={\bf X}(t,\cdot)\# da$ as in (\ref{2defrho}). Assume that $(\bx,\dt\bx) \in \Linf(I\times\Omega)$, with 
$R=\|\bx\|_{\Linf(I\times\Omega)}$, and assume that $\rho \in \Linf(I\times \Rd)$.
Take  
\be
{\bf X}(t)=\nabla\bfi(t)\circ \bg(t),\hspace{1cm} \bg(t)=\nabla\bpsi(t)\circ {\bf X}(t)
\en
to be the polar factorization of ${\bf X}$ as in (\ref{2polar}) where we impose
(\ref{2jauge}). Then 
\begin{enumerate}
\item for a.e. $t\in I$, $\dt\nabla \bfi(t,\cdot)$ is a bounded measure in $\Omega$ with 
\be
\|\partial_t \nabla\bfi\|_{\Linf(I; {\cal M}(\Omega))}
 \leq C(R,d,\Omega)\|\rho\|^{\demi}_{\Linf(I\times B_R)}\|\dt {\bf X}\|_{\Linf(I\times B_r)}
\en
and $\partial_t\bfi \in \Linf(I,L^{1*}(\Omega))$ with $1*=d/(d-1)$.

\item $\bfi$ (resp. $\bpsi$) belongs to $C^{\alpha}(I; C^0(\bar \Omega))$ (resp. to $C^{\alpha}(I; C^0(\bar B_R))$) 
for some $\alpha \in ]0,1[$.

\item  For a.e. $t\in I$, $\partial_t \bg$ is a bounded measure on $\Omega$ with
\be
\|\partial_t \bg\|_{\Linf(I; {\cal M}(\Omega))}\leq C(R,d,\Omega)\|\rho\|_{\Linf(I\times B_R)}
\|\partial_t {\bf X}\|_{\Linf(I\times \Omega)}.
\en
\item If $\rho$ is supported in $\bar\Omega '$ for some open set $\Omega '$, and 
$0< \lambda\leq \rho(\cdot,\cdot)\leq \Lambda$ on $\Omega'$, for some $(\lambda, \Lambda) \in \R^*_+$, 
then there exists $\beta\in ]0,1[$ such that for any $\omega' \subset\subset \Omega'$,
\be
\nabla\bpsi \in C^{\beta}(I;C^0(\omega')),
\en
with $\beta$ depending on $\Lambda/\lambda$.
\item If in addition $\Omega'$ is convex, then there exists $\beta' \in ]0,1[$ such that for any $\omega \subset\subset \Omega$,
\be
\nabla\bfi \in C^{\beta'}(I;C^0(\omega)).
\en
\end{enumerate}
\end{theo}
\vspace{1cm}

\begin{theo}\label{2main2}
Under the assumptions of Theorem \ref{2main1}, 
and assuming that $\rho$ is supported in $\bar \Omega'$, for some open set $\Omega'$, we have:
\begin{enumerate}
\item There exists $\epsilon_0>0$ such that if $|\rho-1|\leq \epsilon < \epsilon_0$ in  $\Omega'$, then there exists $\alpha>0$ (depending on $\epsilon$) such that, for any 
 $w'\subset\subset \Omega'$,
\be
\dt \bpsi \in \Linf(I; C^{\alpha}(\omega')).
\en
If in addition $\Omega'$ is convex,
for any $w\subset\subset \Omega$,
\be
\dt \bfi \in \Linf(I; C^{\alpha}(\omega)).
\en
\item For any $p<2$, there exists $\epsilon(p)>0$ such that, if $|\rho-1|\leq \epsilon(p)$ in $\Omega'$, for any $w'\subset\subset \Omega'$,
\be
\dt \nabla\bpsi \in \Linf(I; L^p(\omega')).
\en
If in addition $\Omega'$ is convex,
for any $w\subset\subset \Omega$,
\be
\dt \nabla\bfi \in \Linf(I; L^p(\omega)).
\en
\end{enumerate}

\end{theo}
\vspace{1cm}
{\it Remark:} The Theorem remains true if one replaces 
the condition $|\rho-1|\leq \epsilon$ by $|\rho-f|\leq \epsilon$
with $f$ a positive continuous function  and the bounds will then depend on the modulus of continuity of $f$ (see \cite{Ca1}) .
\vspace{1cm}

We also state the result in the periodic case:
In this setting we have the following theorem, which is just an adaptation
of the two previous:
\begin{theo}\label{2main2per}
Under the assumptions that $\rho\in \Linf(I\times\Td)$, $\dt\bx \in \Linf(I\times\Td)$, we have:
\begin{enumerate}
\item With the same bounds as in Theorem \ref{2main1}, 
\be
&&\dt \nabla\bfi \in \Linf(I; {\cal M}(\Td)),\\
&& \dt \bg \in \Linf(I; {\cal M}(\Td)),
\en
and for some $\alpha>0$, we have 
\be
\bfi,\bpsi \in C^{\alpha}(I; C^0(\Td)).
\en 
\item If for all $(t,x)\in (I\times\Td)$ we have $0<\lambda\leq \rho(t,x)\leq \Lambda$, then for some $\beta>0$ depending on $(\lambda, \Lambda) \in \R_+^*$, 
\be
\bg, \nabla\bfi, \nabla\bpsi \in C^{\beta}(I; \Linf(\Td)).
\en 
\item There exists $\epsilon_0$ such that if $|\rho-1|\leq \epsilon \leq \epsilon_0$, then for some $\alpha>0$ depending on $\epsilon$,
\be
&&\dt \bpsi \in \Linf(I; C^{\alpha}(\Td)),\\
&&\dt \bfi \in \Linf(I; C^{\alpha}(\Td)).
\en
\item For any $p<2$ there exists $\epsilon(p)$ such that if $|\rho -1| \leq \epsilon(p)$ then 
\be
&&\dt \nabla\bpsi \in \Linf(I; L^p(\Td)),\\
&&\dt \nabla\bfi \in \Linf(I; L^p(\Td)),\\
&&\dt \bg \in \Linf(I; L^p(\Td)).
\en

\end{enumerate}

\end{theo}
{\it Remark:} in this case, the absence of boundary allows to have a bound over $\Td$ and not only interior estimates as in the previous results.


\vspace{1cm}

\subsection{Related results}
\subsubsection*{The linearized Monge-Amp\`ere equation}
The linearized Monge-Amp\`ere  equation $(LMA)$ is a well known equation,  since it is used to carry out the continuity method, in order to obtain classical solutions of the Monge-Amp\`ere equation (see \cite{GT}, chapter 17). 
However for this purpose this is always made in the case where the densities and the domains considered are smooth, and thus the $LMA$ equation is uniformly elliptic.

In the non-smooth case, \cite{CaGu} proved Harnack inequality for solutions of   
\be
M_{ij}\partial_{ij}u=0
\en 
with $M$ the co-matrix of $D^2\Psi$, for some $\Psi$ convex, under the assumption that the measure $\rho=\det D^2\Psi$ satisfies the following absolute continuity condition:
\\
{\it C: For any $0<\delta_1 <1$ there exists $0<\delta_2<1$ such that for any section $S$ and any measurable set $E\subset S$},
\beq
\textrm{if}  \ \frac{|E|}{|S|}\leq \delta_2 \,\textrm{ then } \, \frac{\rho(E)}{\rho{(S)}}\leq \delta_1,
\label{2doubling}
\enq
(a section is a set of the form 
$$
S_t(x_0)=\{ x| \Psi(x)-\Psi(x_0)\leq p\cdot (x-x_0)+t , \  p\in \partial \Psi(x_0) \}).
$$
They showed that the solution of $(\det D^2\Psi) (D^2\Psi)^{-1}_{ij}D_{ij}u=0$ satisfies a Harnack inequality on the sections of 
$\Psi$ 
and subsequently is $C^{\alpha}$. The precise result is stated below (Theorem \ref{2theocagu}).
We will use this result to obtain the first part of Theorem \ref{2main2}.
Note that the condition (\ref{2doubling}) implies $C^{1,\alpha}$ regularity
of the Aleksandrov solution of $\det D^2\Psi=\rho$ (\cite{Ca2}). Note also that the condition (\ref{2doubling}) is satisfied
when the density $\rho$ is bounded between two positive constants.
We will also obtain some results (Theorem \ref{2main1})
in the degenerate case when the condition (\ref{2doubling}) is not satisfied and show in some counterexamples (section \ref{contrex}) that when this condition is not fulfilled, the result of Theorem \ref{2main2} does not hold.

\subsubsection*{Maximum principles for degenerate elliptic equations}
We will use a local maximum principle for degenerate elliptic equations to obtain H\"older continuity in Theorem \ref{2main2}. Consider the problem 
\be
\nabla\cdot (M(x) \nabla u(x)) = \nabla\cdot f(x)
\en
where $M(x)=M_{ij}(x), (i,j)\in [1..d]$ is a symmetric positive semi-definite,
matrix, $f(x)=(f_i(x)) i\in [1..d]$. In the cases we will study, we will not have the usual uniform ellipticity condition 
\be
\lambda I \leq M \leq \Lambda I
\en 
with $I$ the $d\times d$ identity matrix, and  for some positive numbers $\lambda, \Lambda$, but a condition of the form 
\beq\label{2relat}
\lambda(x)I \leq M \leq \Lambda(x)I
\enq
for some non negative measurable functions $\lambda(x), \Lambda(x)$.
Under the assumption that $(\lambda^{-1}(x), \Lambda(x)) \in L^p_{loc}(\Omega)$ for some $p>d$ and that $f\in \Linf$, we can obtain a  bound on the solution $u$ in $\Linf_{loc}$. Properly localized, this bound with the Harnack inequality (Theorem \ref{2theocagu}) will yield H\"older continuity of the solution of the $LMA$ equation (\ref{2diverg}). This type of maximum principles have been already obtained in
\cite{MuSt}, \cite{Tr}, (see also \cite{MuStEr}), and we will use them under the forms of  Theorems \ref{2main5}, \ref{2main4}, and Corollary \ref{2corlocal}. Note however that the condition (\ref{2relat}) is not know by itself to guaranty H\"older continuity of the solution, but only a $\Linf$ bound.

It can be interesting to point out that we will thus use both the divergence and non-divergence structure of the LMA to obtain our results.


\section{Some preliminary results}
In this section we state the results that we are going to need for the proofs of the theorems.
The reader may skip this section and come back to it whenever needed.
Note that all these results can be extended to the periodic case.
\subsection{Regularity for solutions of Monge-Amp\`ere equation}
\begin{theo}\label{2regulier}
Let $\Omega,\Omega'$ be bounded, $C^{\infty}$, strictly convex, and $|\Omega|=1$.
Let $\rho$ be a probability measure in $\bar \Omega'$, belong to $C^{\infty}(\bar\Omega')$, and satisfy
$0 < \lambda \leq \rho(t,x) \leq \Lambda$ for some pair $(\lambda, \Lambda)$. 
Then there exists a unique (up to a constant) solution of  
\begin{eqnarray*} 
&&\det D^2\bpsi=\rho \label{2maper1},\\
&&\nabla\bpsi \textrm{ maps } \Omega' \textrm{ to }  \Omega \label{2maper2}, 
\end{eqnarray*} 
in the sense of (\ref{2defPsi}).
The solution $\bpsi$ belongs to $C^{\infty}(\bar \Omega')$, and
$\bfi$, defined as in (\ref{2defPhi}), belongs to $C^{\infty}(\bar \Omega)$.
\end{theo} 
For this the reader can refer to \cite{Ca1}-\cite{Ca3}, \cite{De}, \cite{U}.

\vspace{1cm}
\noindent
The next Theorem can be found in  \cite{Ca2}, \cite{Ca3}, \cite{Ca4}.
\begin{theo}\label{2caf}
Let $\rho$ be supported in $\bar\Omega'$ with $\Omega'$ open, satisfy $0<\lambda \leq\rho\leq \Lambda$, and let $\bpsi$ be solution of  
\be
&&\det D^2 \bpsi =\rho,\\
&&\nabla\bpsi \textrm{ maps } \Omega' \textrm{ to }  \Omega,
\en in the sense of (\ref{2defPsi}) with $\Omega$ convex.
Then for some $\alpha \in ]0,1[$ depending on $\Lambda/\lambda$, $\bpsi \in C^{1,\alpha}_{loc}(\Omega')$. If moreover $\Omega'$ is also convex then
$\bpsi$ (resp.  its Legendre transform $\bfi$) is in $C^{1,\alpha}(\bar\Omega')$ (resp. in $C^{1,\alpha}(\bar\Omega)$).
\end{theo}

\vspace{1cm}
\noindent
The next Theorem can be found in \cite{Ca1}.
\begin{theo}\label{2W2p}
Let $\Omega$ be normalized so that $B_1 \subset \Omega \subset B_d$.
Let $\bpsi$ be a convex Aleksandrov solution of
\be
&&\det D^2\bpsi =\rho,\\
&& \bpsi = 0 \text{ on } \partial\Omega.
\en
Then for every $p< \infty$ there exists $\epsilon(p)$ such that
if $|\rho-1|\leq \epsilon(p)$ then $\bpsi \in W^{2,p}_{loc}(\Omega)$ and  
\be
\|\bpsi\|_{W^{2,p}(B_{1/2})}\leq C(\epsilon).
\en
\end{theo}
{\it Remark 1:} This implies also, maybe for a smaller value of $\epsilon(p)$ that one can also have
$\ds\|D^2\bpsi^{-1}\|_{L^p(B_{1/2})} \leq C'(\epsilon)$. 
\\
{\it Remark 2:} The theorem remains true if one replaces $|\rho-1|\leq \epsilon$ by $|\rho-f|\leq \epsilon$, for some
continuous positive $f$, 
and the bounds depends on the modulus of continuity of $f$.


\vspace{1cm}

\subsection{The linearized Monge-Amp\`ere equation}
We state here the result of \cite{CaGu} evoked in the previous section: 
\begin{theo}\label{2theocagu}
Let $\Omega$ be a domain in $\Rd$, let $U$ be an Aleksandrov  solution in $\Omega$ of 
\be
\det D^2 U = \mu
\en
where $\mu$ the  satisfies the condition (\ref{2doubling}). 
Let $w$ be a solution in $\Omega$ of the linearized homogeneous Monge-Amp\`ere equation
\be
A_{ij} \partial_{ij}w=0
\en
where $A_{ij}$ is the co-matrix of $D^2 U$,
let $R>0$ and $y\in \Omega$ be such that $B_R(y) \subset \Omega$,
 then for some $\beta<1$ depending only on the condition (\ref{2doubling}), for any $r<R/4$, 
\be
osc(r/2)\leq \beta osc(r),
\en
where
\be
&& osc(r)=M(r)-m(r),\\
&& M(r)=\sup_{B_r(y)} w, \hspace{1cm} m(r)=\inf_{B_r(y)} w.
\en
\end{theo}


\subsection{Maximum principle for degenerate elliptic equations}
We give here some results concerning degenerate elliptic equations of the form
\beq\label{2divforme}
\nabla \cdot (M(x)\nabla u(x))= \nabla\cdot f(x)
\enq
where $M$ is  symmetric non-negative matrix, $f=(f_i), i=1..d$. The equation can be written  
$\partial_i(M_{ij}\partial_j u)=\partial_i f_i$ with summation over repeated indices.
The usual strict ellipticity condition
\be
\lambda |\xi|^2 \leq M_{ij}\xi_i\xi_j \leq \Lambda |\xi|^2 \text{ for all } \xi \in \R^d,
\en
is replaced by the following
\be
\sum_{i,j=1}^d |M_{ij}| + |M^{ij}| \in L^p_{loc}(\Omega)  \text{ for some } p,
\en
where $M^{ij}$ denotes the inverse matrix of $M$. This is equivalent to the condition that there exists $\lambda(x), \Lambda(x)$ such that $\lambda^{-1}, \Lambda$ are in $L^p_{loc}(\Omega)$ and such that $\lambda(x) I \leq M(x) \leq \Lambda(x) I$, in the sense of symmetric matrices.

The class of admissible test functions
is 
\be
{\cal C}(\Omega)= \{ v \in W^{1,1}_0(\Omega), \ M^{1/2} \nabla v \in L^2 (\Omega)\}.
\en
A subsolution (resp. supersolution) $u$ of (\ref{2divforme}) is defined by the condition that for all non-negative $v\in {\cal C}(\Omega)$, 
\be
\int_{\Omega} \nabla v M \nabla u - \nabla v \cdot f  \leq (\geq) 0.  
\en
Then, following \cite{MuSt} and \cite{Tr}, we have the following results:

\subsubsection*{Bound for Dirichlet boundary data}

We denote by $S_d^+$ the set of $d\times d$ non negative symmetric matrices.
\begin{theo}\label{2main5}
Let $M: \Omega \to S_d^+$ be such that $M^{-1}$ is in $L^p(\Omega; S_d^+)$ for some $p>d$. Let $f$ be in $\Linf(\Omega;\Rd)$. Let $u$ be a subsolution (supersolution) of 
\be
\nabla \cdot (M(x)\nabla u(x))= \nabla\cdot f(x)
\en
in $\Omega$,
satisfying $u\leq 0$ ($u\geq 0$) on $\partial\Omega$. Then 
\be
\sup_{\Omega} u (-u) \leq C(\|u^+(u^-)\|_{L^{a_0}(\Omega)} + \|f\|_{\Linf(\Omega)})
\en
where $C,C$ depends on $|\Omega|, a_0>0, p>d, \|M^{-1}\|_{L^p(\Omega)}$. 
\end{theo}

\vspace{1cm}
\noindent
This maximum principle can be precised in the following corollary, that will be crucial for the proof of H\"older continuity in Theorem \ref{2main2}.  
\begin{cor}\label{2corlocal}
Under the previous assumptions,  
for $y\in\Omega$, $B_R(y) \subset \Omega$,
if $u$ is a subsolution (supersolution) in $B_R$ of (\ref{2divforme}) and $u\leq 0$ ($u\geq 0$) on $\partial B_R$, then
\be
\sup_{B_R} u (-u) \leq C\|M^{-1}\|_{L^p(B_R)}\|f\|_{\Linf(B_R)} R^{\delta},
\en
where  $\ds\delta=1-\frac{n}{p}$.
\end{cor}

\subsubsection*{Bound without boundary data}
Here we state a maximum principle that does not depend on the boundary
data. Note that here we need to control the norm of both $M$ and $M^{-1}$ whereas we only needed to control $M^{-1}$ above.
\begin{theo}\label{2main4}
Let $M: \Omega \to S_d^+$ be such that $M, M^{-1}$ are both in $L^p_{loc}(\Omega)$, with $p>d$. Let $f$ be in $\Linf(\Omega)$. Let $u$ be a subsolution of 
\be
\nabla \cdot (M(x)\nabla u(x))= \nabla\cdot(f(x))
\en
in $\Omega$. 
Then we have for any ball $B_{2R} \subset \subset \Omega$ and $a_0>0$
\be
\sup_{B_R(y)} u \leq C_1\|u^+\|_{L^{a_0}(B_{2R}(y))} + C_2 k
\en
where $k= \|f\|_{\Linf(B_{2R})}$, $C_1, C_2$ depend on $R, a_0, p, \|M\|_{L^p(B_{2R})},\| M^{-1}\|_{L^p(B_{2R})}$. 
\end{theo}

\subsection{Convex functions and Legendre transforms}
We state first the following classical lemma on convex functions:
\begin{lemme}\label{2D2phiL1}
Let $\varphi$ be a convex function from $\Rd$ to $\R$, globally Lipschitz with Lipschitz constant
$L$. Then we have
\be
\|D^2\varphi\|_{{\cal M}(B_R)}\leq C(d)R^{d-1} L.
\en
\end{lemme}
{\it Proof:} 
we have 
\be
\|D^2\varphi\|_{{\cal M}(B_R)}&\leq &C\int_{B_R}\Delta\varphi \\
&=& \int_{\partial B_R}\nabla\varphi\cdot n \\
&\leq& C(d) R^{d-1}L.
\en

$\hfill \Box$

\vspace{1cm}
\noindent
We recall here some useful properties of the Legendre transform.
Let $\Omega$ be a convex domain, let $\phi:\Omega \mapsto \R$ be $C^1$ convex.
Let $\phi^*$ be its Legendre transform 
defined by 
\be
\phi^*(y)=\sup_{x\in \Omega} x\cdot y - \phi(x).
\en
Then, for all $x\in \Omega$,   
\be
&&\nabla\phi^*(\nabla\phi(x))=x.
\en
If moreover $\phi$ is $C^2$ strictly convex, then, for all $x\in \Omega$,
\beq\label{2legendreinverse}
D^2\phi^*(\nabla\phi(x))=D^2\phi^{-1}(x).
\enq
From this we deduce the following lemma: 
\begin{lemme}
Let $\Omega$ be convex, let $(t,x) \mapsto \bfi(t,x): I\times\Omega \mapsto \R$ and 
$(t,y) \mapsto \bpsi(t,y): I\times\R^d \mapsto \R$ be such that
\begin{enumerate}
\item $\nabla\bfi$  (resp. $\nabla\bpsi$) belongs to $C^1(I\times \Omega)$ (resp. belongs to $C^1(I\times \R^d)$),
\item for all $t\in I$, $\bfi(t,\cdot)$ is convex and $\bpsi(t,\cdot)$ is the Legendre transform of $\bfi(t,\cdot)$.
\end{enumerate}
then for every $(t,x) \in I\times\Omega$,
\beq
&&\bfi(t,x)+\bpsi(t,\nabla\bfi(t,x))=x\cdot\nabla\bfi(t,x)\label{2identite0},\\
&&\partial_t\bfi+\partial_t\bpsi(\nabla\bfi)=0\label{2identite1},\\
&&\partial_t\nabla\bfi+ D^2\bfi\partial_t\nabla\bpsi(\nabla\bfi)=0\label{2identite2}.
\enq
\end{lemme}
{\it Proof:} 
the first identity expresses just the fact that $\bfi(t,\cdot), \bpsi(t,\cdot)$ are Legendre transforms of each other (see (\ref{2deflegendre})), then the two other come by differentiating with respect to time and then to space.

$\hfill \Box$


\section{Approximation by smooth functions}
\subsection{Construction of smooth solutions.} 
In this section we build an adequate smooth approximation of the problem.
More precisely, given a mapping $\bx(t)$ and $\rho(t)=\bx(t)\#da$,  we construct an associated pair $(\rho, v)$ satisfying 
\beq\label{2continuite}
\dt\rho + \nabla\cdot(\rho v)= 0
\enq
and  then  find a ``good'' regularization of $(\rho,v)$. 
One of the problems is the following: it is known from a counterexample by Caffarelli (see \cite{Ca3}), that
when transporting a (smooth) density $\rho_1$ onto another (smooth) density $\rho_2$ by the gradient of a convex function, 
one can not expect the convex function to be $C^1$  
unless $\rho_2$ is supported and positive in a convex set. Therefore it is not enough to only regularize 
(by convolution for example) the density $\rho=\bx{\#}da$, we must also approximate it by a density supported in a
convex set.
\\
The density $\rho$ and $\dt \rho $ are  constructed from ${\bf X}, \dt {\bf X}$  respectively 
by the following procedure:
\be
\forall f \in C^1_b(\Rd),  &&\int_{\Rd} \rho(t,x)f(x) dx =\int_{\Omega} f({\bf X}(t,a))\, da\\
&&\int_{\Rd} \dt \rho(t,x)f(x) dx=\int_{\Omega} \nabla f({\bf X}(t,a))\cdot \dt {\bf X}(t,a)\, da.
\en
To define $v$ such that $\dt\rho + \nabla\cdot (\rho v)=0$,  
we define the product $\rho v$ as follows: 
\be
\forall \phi \in C^0_b(I\times\Rd; \R^d),\int_{I\times \Rd}\rho v\cdot \phi \,  dt dx = \int_{I\times \Omega}\phi({\bf X}(t,a))\cdot \dt {\bf X}(t,a) \, dt da. \en
Since $\dt X\in \Linf$, $v$ is well defined $d\rho$ a.e. and we have
\be
\|v(t,\cdot)\|_{\Linf(\Rd, d\rho(t))}\leq \|\dt {\bf X}(t,.)\|_{\Linf(\Omega)}.
\en
Now we construct $(\rho_n,v_n)$ a smooth approximating sequence for
$(\rho,v)$ as follows: (remember that we have taken $\rho(t,\cdot)$ to be supported in $B_R$ at any time $t\in I$). We take $\eta \in C^{\infty}_c$ a standard  convolution kernel, of integral 1, supported in $B(0,1)$ and positive. Take $\eta_n=n^d\eta(n x)$. We also note $\chi_{R+1/n}$ the characteristic function of the ball $B(0, R+1/n)$. Let
\be
&&\rho_n=(\frac{1}{n}\chi_{R+1/n}+\eta_n*\rho)c_n,\\
&&v_n=c_n\frac{\eta_n*(\rho v)}{\rho_n}, 
\en 
with $c_n$ chosen such that $\rho_n$ remains a probability measure. (Note that $c_n$ is close to 1 for $n$ large). The purpose of this construction is to have the following properties:
\begin{enumerate}
\item $\|\rho_n,v_n\|_{\Linf}\leq\|\rho,v\|_{\Linf}$,  
\item $\rho_n, v_n$ satisfy the continuity equation (\ref{2continuite}),
\item  $\rho_n$ is supported and strictly positive in $B(0, R+1/n)$, and belongs to $C^{\infty}(\bar B(0, R+1/n))$ .
\item If $\bfi_n(t), \bpsi_n(t)$ are associated to $\rho_n(t)$ through (\ref{2defPhi},\ref{2defPsi}), then, for every $t\in I$,  $\bfi_n(t)$ converges uniformly on compact sets of $\Omega$ to $\bfi(t)$
and $\bpsi_n(t)$ converges uniformly on compact sets of $\Rd$ to $\bpsi(t)$. This last result can be found in \cite{Br1}. Therefore, $\dt\bfi_n, \dt\bpsi_n$ will converge in the distribution sense to $\dt\bfi, \dt\bpsi$.
\end{enumerate}

\vspace{1cm}
\noindent
Now we have the following regularity result, for smooth densities. Note that this result
will only be used to legitimate the forthcoming computations, and not as an a-priori bound.  

\begin{prop} \label{2main3a}
let $I, \Omega$ be as above, let $\Omega'$ be $C^{\infty}$ strictly convex. 
For any $t\in I$, let $\rho(t,\cdot)$ be a probability density in  $\bar\Omega'$, strictly positive in $\bar\Omega'$ 
with $\rho \in C^{\infty}(I\times \bar\Omega')$. 
Let, for all $t$,  $\bfi(t,\cdot), \bpsi(t,\cdot)$ be 
as in (\ref{2defPhi},\ref{2defPsi}) with $(\rho(t), \Omega)$. Then, for any $0<\alpha<1$,
\be
\dt\bfi \in \Linf(I,C^{2,\alpha}(\bar\Omega)),  \hspace{1cm} \dt\bpsi \in\Linf(I,C^{2,\alpha}(\bar \Omega')).
\en
\end{prop}

\vspace{1cm}

\noindent
{\it Proof of Proposition \ref{2main3a}:}
Theorem \ref{2regulier} implies that for all $t$, $D^2\bpsi$ (resp. $D^2\bfi$)
belongs to $C^{\infty}(\bar \Omega')$ (resp. belongs to $C^{\infty}(\bar\Omega)$).
\\
Now we wish to solve $\det D^2 \bpsi(t)=\rho(t)$ with $t$ near $t_0$.
We write a priori  $\bpsi(t)=\bpsi(t_0) + (t-t_0)u + o(|t-t_0|)$, 
for some $u$, then we have
\be
\det D^2 \bpsi(t)=\det D^2\bpsi(t_0) + (t-t_0)\text{trace}(M D^2 u) + o(|t-t_0|)
\en
where $M$ is the comatrix of $D^2\bpsi$  defined by 
\be
M(t,x)=\det D^2 \bpsi(t,x) \left( D^2\bpsi(t,x)\right)^{-1}.
\en
Note that $M$ belongs to $C^{\infty}(\bar \Omega')$ and is uniformly elliptic.
Let us now show that $\dt \bpsi$ can indeed be sought as the solution of
\be
\text{trace}(M D^2 u)=\dt\rho
\en 
with a suitable boundary condition. 
For this we introduce $h$ a defining function for $\Omega$, ({\it i.e.} $h\in C^{\infty}(\bar\Omega)$ is strictly convex and vanishes on $\partial\bar\Omega$, we can also impose $|\nabla h|_{\partial\bar\Omega}\equiv 1$). 
The condition $\nabla\bpsi$ maps $\Omega'$ on $\Omega$ can be replaced by $h(\nabla\bpsi)=0$ on $\partial \Omega'$.  Now consider the operator
\be
{\cal F}: \psi \mapsto \left(\det D^2\psi, h(\nabla\psi)|_{\partial\Omega'}\right)
\en defined on $\{\psi\in C^{2,\alpha}(\bar \Omega'), \psi \text{ convex }\}$ and ranging in $C^\alpha(\bar \Omega')\times C^{1,\alpha}(\partial \Omega')$.
First note that a smooth solution of 
\beq\label{2ilefo}
{\cal F}(\psi)=(\rho(t), 0)
\enq
will satisfy (\ref{2defPsi}) and thus co\"incide (up to a constant) with $\bpsi(t)$. 
We now solve (\ref{2ilefo}) around $t_0$ by the implicit function Theorem. 
The derivative of ${\cal F}$ at $\bpsi$ is defined by
\be
d{\cal F}(\bpsi) u = (I(u), B(u))=\left(M_{ij} \partial_{ij}u, h_i(\nabla\bpsi) \partial_{i}u  \right).
\en
The operator $I=M_{ij}\partial_{ij}$ is uniformly elliptic with coefficients $M_{ij}$ in $C^{\infty}(\bar \Omega')$.
We need also to show that the boundary operator $B$ is strictly oblique:
First, note that $\nabla h=\vec{n}_1$ on $\partial \Omega$, where  $\vec{n}_1$ is the outer unit normal to $\partial \Omega$. 
Moreover, if  $\vec{n}_2$ is the outer unit normal to $\partial \Omega'$, 
it has been established in \cite{Ca3}, \cite{De}, \cite{U}, 
that there exists a constant $C$ depending on $\Omega, \|\rho\|_{C^2(\bar \Omega')}$, and therefore uniform on $I$, such that 
\be
\vec{n}_2\cdot \vec{n}_1(\nabla\bpsi) \geq C>0.
\en
Thus the boundary condition is strictly oblique, uniformly with respect to $t$.
It has been established in \cite{De}, p. 448, that the equation 
\be
d{\cal F}(\bpsi) u= (\mu, 0)
\en 
with $\mu \in C^{\alpha}(\bar \Omega')$
is solvable up to an additive constant if $\int_{\Omega'} \mu = 0$. 
This condition is met by $\dt \rho$, since $\int \rho(t,x) \ dx\equiv 1$.
\\
We conclude that the operator $d{\cal F}(\bpsi)$ is invertible on 
the set $$\left\{\{\mu \in C^{\alpha}(\bar \Omega'), \int \mu =  0\} \times \{\nu=0\}\right\}$$ {\it i.e.} for each $\mu\in C^{\alpha}(\bar\Omega')$, with $\int_{\Omega'} \mu =0$, there exists a unique up to a constant solution 
$u$ of $d{\cal F}(\nabla\bpsi)u= (\mu ,0)$. Moreover, following \cite{GT}, Theorem 6.30, $u$ belongs to $C^{2,\alpha}(\bar \Omega')$. 
\\
Therefore we can apply the implicit function Theorem and solve 
${\cal F}(\bpsi(t))=(\rho(t),0)$ for $t$ near $t_0$.
By uniqueness of the solution of (\ref{2defPsi}), this solution will co\"incide with the solution of Theorem \ref{2regulier}.
As we have built it, $\dt\bpsi(t,\cdot)=u$ is the unique (up to a constant) 
solution of 
\begin{eqnarray} 
&&\text{trace}\left(M D^2 u\right)=\partial_t\rho\,\,\textrm{ in } \Omega',\\ 
&&\nabla u \cdot \vec{n_1}(\nabla\bpsi)=0\,\,\text{ in }\partial \Omega',     
\end{eqnarray}
and  since $\dt\rho\in C^{\infty}(\bar\Omega')$, $\dt\bpsi$ belongs to $C^{2,\alpha}(\bar \Omega')$ for any $\alpha<1$.
\\
We also have, using the identity (\ref{2identite1})
\be
\dt \bfi + \dt \bpsi(\nabla\bfi)=0.
\en
therefore
$\dt\bfi \in C^{2,\alpha}(\bar\Omega)$ for any $\alpha<1$.
\\ 
This achieves the proof of Proposition \ref{2main3a}. 

$\hfill\Box$

\section{Proof of Theorem \ref{2main1}}

Theorem \ref{2main1} will be deduced through approximation from the following proposition: 


\begin{prop}\label{2main3b}
Let $\rho$ satisfy the assumptions of Proposition \ref{2main3a} above, with $\Omega'=B_R$,  and $\bfi, \bpsi$ be as in (\ref{2defPhi}, \ref{2defPsi}).
Let  $v(t,x)\in \R^d$  be a smooth vector field on $\bar B_R$ and satisfy on $I\times \bar B_R$ 
\beq 
\dt \rho + \nabla\cdot(\rho v)=0\label{2conti}. 
\enq 
Take $1\leq p,r \leq \infty$, $\frac{1}{r}+\frac{1}{r'}=1$, $q=\frac{2p}{1+p}$. Then for any $t\in I$, for any $\omega \subset \Omega$ we have:
\beq\label{2dtgradphi}
\|\partial_t\nabla\bfi\|_{L^q(\omega)}\leq \left(\|\rho |v|^2\|_{L^{r'}}\|D^2\bpsi\|_{L^r(B_R)} \|D^2\bfi\|_{L^p(\omega)} \right)^{1/2},
\enq
\\
which implies in particular 
\beq\label{2dtPhiL1}
\|\dt\nabla\bfi\|_{L^1(\Omega)}\leq C(R,d,\Omega)\left(\|\rho |v|^2\|_{\Linf}\right)^{1/2}, 
\enq
and for any $t\in I$, for any $\omega' \subset B_R$ we have: 
\beq\label{2dtgradpsi}
\hspace{.5cm}\left[\int_{\omega'} \rho |\partial_t\nabla\bpsi|^q\right]^{1/q} \leq  \left(\|\rho |v|^2\|_{L^{r'}}\|D^2\bpsi\|_{L^r(B_R)} 
 \left[\int_{\omega'} \rho |D^2\bpsi|^{p}\right]^{1/p}\right)^{1/2},
 \enq
\\
which  implies in particular
\beq\label{2dtPsiL1}
\int_{\Rd} \rho|\dt\nabla\bpsi| \leq C(R,d,\Omega)\|\rho\|^{\demi}_{\Linf(\Rd)}\|\rho |v|^2\|^{\demi}_{\Linf(\Rd)}.
\enq
\end{prop}

\vspace{1cm}

\noindent
{\it Proof of Proposition \ref{2main3b}:}
\\
Using Proposition \ref{2main3a}, we can perform the following computations. 
We have from (\ref{2defPhi})
\be
\int_{\Rd} \dt \bpsi \rho = \int_{\Omega}\dt \bpsi (\nabla\bfi)
\en
Then we use the continuity equation:
\be
\dt\rho + \nabla\cdot(\rho v)=0
\en
which implies for any smooth $f$ 
\be
\int_{\Rd} f \dt \rho = \int_{\Rd} \rho v\cdot \nabla f.
\en
We obtain 
\be
\int_{\Rd} \partial_t\bpsi \dt\rho &=& \int_{\Rd}\partial_t\nabla\bpsi\cdot\rho
v\\ 
&=&\int_{\Omega}\partial_t\nabla\bpsi(\nabla\bfi)\cdot\dt\nabla\bfi\\
&=&-\int_{\Omega}\partial_t\nabla^t\bpsi(\nabla\bfi)\cdot D^2\bfi\cdot\partial_t \nabla\bpsi
(\nabla\bfi) 
\en 
where we have used (\ref{2identite2}).
Since we can write $\sqrt{D^2\bfi}$ because this is a
positive symmetric matrix, we have 
\be
\|\sqrt{D^2\bfi} \ \partial_t\nabla\bpsi(\nabla\bfi)\|_{L^2(\Omega)}^2
&=&-\int_{\Rd} \rho \partial_t\nabla\bpsi\cdot  v\\ 
&=&-\int_{\Omega}\partial_t\nabla\bpsi(\nabla\bfi)\cdot v(\nabla\bfi)\\ 
&=&-\int_{\Omega}\sqrt{D^2\bfi}\partial_t\nabla\bpsi(\nabla\bfi) \cdot \sqrt{D^2\bfi}^{-1} v(\nabla\bfi).
\en 
This implies that 
\beq\label{2rst}
\|\sqrt{D^2\bfi} \ \partial_t\nabla\bpsi(\nabla\bfi)\|_{L^2(\Omega)}
\leq \|\sqrt{D^2\bfi}^{-1}v(\nabla\bfi)\|_{L^2(\Omega)}.
\enq

In order to estimate the right hand side, we write 
\beq 
\|\sqrt{D^2\bfi}^{-1}v(\nabla\bfi)\|_{L^2(\Omega)}&=&\left(\int_{\Omega}
v^t(\nabla\bfi)\cdot(D^2\bfi)^{-1}\cdot v(\nabla\bfi)\right)^{1/2}\nonumber\\
& =&\left(\int_{\Rd}
\rho v^t\cdot (D^2\bfi(\nabla\bpsi))^{-1}\cdot v\right)^{1/2}\nonumber\\ &=&\left(\int_{\Rd}
\rho v^t\cdot D^2\bpsi \cdot v\right)^{1/2}\nonumber\\
&\leq &
\left(\|D^2\bpsi\|_{L^r(B_R)}\|\rho v^2\|_{L^{r'}(B_R)}\right)^{1/2}.\label{2above} 
\enq 
In the second line we have used  $D^2\bfi(\nabla\bpsi)=(D^2\bpsi)^{-1}$.
From (\ref{2identite2}), 
\be
\|\sqrt{D^2\bfi}^{-1}\partial_t\nabla\bfi\|_{L^2(\Omega)}&=&\|\sqrt{D^2\bfi}\partial_t\nabla\bpsi(\nabla\bfi)\|_{L^2(\Omega)}\\
&\leq&\|\sqrt{D^2\bfi}^{-1}v(\nabla\bfi)\|_{L^2(\Omega)}. 
\en 
Writing 
$$\partial_t\nabla\bfi=\sqrt{D^2\bfi}^{-1}\sqrt{D^2\bfi}\partial_t\nabla\bfi,$$
and, using H\"older's inequality, we obtain for $\omega\subset\Omega$
\be
\|\partial_t\nabla\bfi\|_{L^q(\omega)}&\leq&\|\sqrt{D^2\bfi}^{-1}\partial_t\nabla\bfi\|_{L^2(\omega)}\|\sqrt{D^2\bfi}\|_{L^s(\omega)}\nonumber\\
&\leq& \left(\|\rho |v|^2\|_{L^r(B_R)}\|D^2\bpsi\|_{L^{r'}(B_R)}\|D^2\bfi\|_{L^{s/2}(\omega)}\right)^{1/2}
\en 
with $q=\frac{2s}{2+s}.$ 
By taking $p:=s/2$
we have 
\be
\|\partial_t\nabla\bfi\|_{L^q(\omega)}&\leq& \left(\|\rho |v|^2\|_{L^r(B_R)}\|D^2\bpsi\|_{L^{r'}(B_R)}\|D^2\bfi\|_{L^{p}(\omega)}\right)^{1/2}
\en
and $q=\frac{2p}{1+p}$. This proves (\ref{2dtgradphi}).
\bigskip
\noindent
To obtain a bound on $\partial_t\bpsi$ we write
\be
\int_{\Rd} \rho \left|\sqrt{D^2\bfi (\nabla\bpsi)}\partial_t\nabla\bpsi\right|^2
&=&\int_{\Rd} \rho \partial_t\nabla^t\bpsi\cdot D^2\bfi (\nabla\bpsi)\cdot\partial_t\nabla\bpsi \\
&=& \int_{\Omega} \partial_t\nabla^t\bpsi(\nabla\bfi)\cdot D^2\bfi \cdot\partial_t\nabla\bpsi(\nabla\bfi)\\ 
&\leq& \|D^2\bpsi\|_{L^r(B_R)}\|\rho |v|^2\|_{L^{r'}(B_R)} 
\en
from (\ref{2rst}) and (\ref{2above}).
Then using H\"older's inequality, with $q=\frac{2s}{2+s}$, we obtain for $\omega' \subset B_R$,
\be
&&\left[\int_{\omega'} \rho |\partial_t\nabla\bpsi|^q\right]^{1/q} \leq  \\
&& \left[\int_{\omega'} \rho \left| \sqrt{D^2\bfi (\nabla\bpsi)}\partial_t\nabla\bpsi\right|^2\right]^{1/2} \left[\int_{\omega'} \rho \left|[D^2\bfi (\nabla\bpsi)]^{-1}\right|^{s/2}
\right]^{1/s}.
\en
The first factor of the right hand product has been estimated above, and the second is equal to
$\ds\left(\int \rho |D^2\bpsi|^{s/2}\right)^{1/s}$.
We conclude that 
\be
\left[\int_{\omega'} \rho |\partial_t\nabla\bpsi|^q\right]^{1/q} \leq 
 \left[ \|D^2\bpsi\|_{L^r(B_R)}\|\rho |v|^2\|_{L^{r'}(B_R)} \right] ^{1/2}
\left[\int_{\omega'} \rho |D^2\bpsi|^{s/2}\right]^{1/s}. 
\en
Taking again $p:=s/2$, we have proved (\ref{2dtgradpsi}).
\\
The bounds (\ref{2dtPhiL1}, \ref{2dtPsiL1}) are obtained as follows: we know from Lemma \ref{2D2phiL1}
that 
\be
&&\|D^2\bpsi\|_{L^1(B_R)}\leq C(R,d,\Omega),\\
&&\|D^2\bfi\|_{L^1(\Omega)}\leq C(R,d,\Omega).
\en
Taking in (\ref{2dtgradphi}, \ref{2dtgradpsi}) $r=+\infty, r'=1, p=1$ we obtain the desired bounds.
This ends the proof of Proposition \ref{2main3b}.

$\hfill\Box$

\subsection{Proof of Theorem \ref{2main1}}

\subsubsection*{Proof of the bound on $\dt\nabla \bfi$}
Here we prove points 1,2,4,5 of Theorem \ref{2main1}.
To obtain point 1, we just need to pass to the limit in the estimate (\ref{2dtPhiL1}).
We need to have $\ds \liminf \|\rho_n |v_n|^2 \|_{\Linf}\leq\|\rho |v|^2 \|_{\Linf}$: to prove this, 
notice that $\ds F(\rho,v)=\rho |v|^2/2=\frac{(\rho |v|)^2}{2\rho}$ is a convex functional in $(\rho v, \rho)$ since it is expressed as:
\be
\frac{(\rho |v|)^2}{2\rho}=\sup_{c+|m|^2/2\leq 0}\{ \rho c + \rho v \cdot m\}.
\en
Then since $\ds \rho_n v_n = c_n \eta_n*(\rho v), \rho_n = c_n(\frac{1}{n}+ \eta_n*\rho)$ we get that
\be
F(\rho_n, \rho_n v_n)\leq c_n \eta_n * F(\rho, \rho v) \leq c_n\|\rho \frac{|v|^2}{2} \|_{\Linf}
\en
and letting $n \to \infty$:
\be
\|\partial_t \nabla\bfi\|_{{\cal M}(\Omega)}
&\leq& \left(\|\rho |v|^2 \|_{\Linf}\right)^{\demi}C(R,d,\Omega)\\
 &\leq& \|\rho\|^{\demi}_{\Linf(B_R)}\|v\|_{\Linf(B_r,d\rho)}C(R,d,\Omega).
\en
Since we impose   
$\ds \int_{\Omega}\bfi(t,x) \ dx \equiv 0$, and since $\Omega$ is convex, (note that since $\ds\dt \bfi_n \notin W^{1,1}_0$, a condition of this type is necessary, see \cite{GT}, chap. 7)
by Sobolev imbeddings  we get 
also a bound on $\|\partial_t\bfi_n\|_{L^{1*}(\Omega)}$.
This proves the first point of Theorem \ref{2main1}.

\vspace{.5cm}

\noindent
Then we obtain 
points 2,4,5 by the following interpolation lemma:
\vspace{.5cm}
\begin{lemme}\label{2Calpha}
Let $\bfi_1$ and $\bfi_2$ be two $R-Lipschitz$ convex functions on $\Omega$ convex.
Then 
\\
1- there exists
$C,\beta>0$ depending on $(\Omega, R, d, p)$ such that
\be
\|\bfi_1-\bfi_2\|_{\Linf(\Omega)}\leq C \|\bfi_1-\bfi_2\|^{\beta}_{L^p(\Omega)}.
\en
\\
2- If moreover $\bfi_1\in C^{1,\alpha}$ for some $0<\alpha<1$ then  there exists 
$C', \beta'>0$ depending also on $\alpha, \|\bfi_1\|_{C^{1,\alpha}}$,  such that, if 
$\Omega_\delta=\{x\in \Omega, d(x,\partial \Omega)\geq \delta\}$, with $\delta$ going to 0 with $\|\bfi_1-\bfi_2\|_{L^p(\Omega)}$,  then
\be
\|\nabla \bfi_1-\nabla \bfi_2\|_{\Linf(\Omega_\delta)}\leq C'\|\bfi_1-\bfi_2\|_{L^p(\Omega)}^{\beta'}.
\en
\end{lemme}

\vspace{1cm}

\noindent
{\it Proof:} Suppose that $\ds\int_{\Omega}|\bfi_1-\bfi_2|^p \leq \epsilon^p $. 
Choose a point inside $\Omega$ (say 0) such that $|\bfi_1(0)-\bfi_2(0)|=M$. $\bfi_1$ and $\bfi_2$ are globally Lipschitz  with Lipschitz constant bounded by $R$.  
On $B_{M/2R}(x) \cap \Omega$ we have 
$|\bfi_1-\bfi_2|(x)\geq M/2$
 and thus 
\be
\int_{B_r}|\bfi_1-\bfi_2|^p \geq \textrm{vol}(\Omega\cap B_{M/2R}(x)) (M/2)^p.
\en
Next note that for $\Omega$ convex, $M$ small enough,  for any $x\in \Omega$, 
$\textrm{vol}(\Omega\cap B_{M/2R}(x))\geq C_{\Omega}\textrm{vol}(B_{M/2R}(x))$. Finally we have 
\beq\label{2M}
\epsilon^p\geq\int_{\Omega}|\bfi_1-\bfi_2|^p \geq C(\Omega,R,d) M^ {p+d},
\enq
and thus 
\be
M\leq C'(\Omega,R,d)\left[\int_{B_r}|\bfi_1-\bfi_2|^p\right]^{\frac{1}{p} \frac{p}{p+d}}, 
\en
which gives the first part of the lemma, with $\ds\beta=\frac{p}{p+d}$.
\\
Now suppose that $|\nabla\bfi_1(0)-\nabla\bfi_2(0)|=M$. One can also set $\bfi_1(0)=0, \nabla\bfi_1(0)=0$. 
We know that $\bfi_1$ is $C^{1,\alpha}$ thus $\bfi_1(x)\leq C|x|^{1+\alpha}$. It follows that going in the direction of $\nabla\bfi_2$ one will have
$$\bfi_2(x)-\bfi_1(x)\geq M|x| - C|x|^{1+\alpha}+\bfi_2(0).$$ 
Keeping in mind that  $|\bfi_1(x)-\bfi_2(x)|\leq C\epsilon^\beta $ yields
$M|x| - C|x|^{1+\alpha}\leq C\epsilon^\beta $. The maximum of the left hand side
is attained  for $|x|=\left(\frac{M}{(1+\alpha)C}\right)^{1/\alpha}$, and is equal to 
$\left(\frac{M}{(1+\alpha)C}\right)^{1/\alpha}\frac{\alpha}{1+\alpha}M$.
Therefore we have 
$$M\leq C\epsilon^{\beta'}$$
in $\Omega_\delta$ with $\delta=\delta(\epsilon)$ going to 0 as $\epsilon$ goes to  0
and with $\beta'=\frac{\alpha \beta}{1+\alpha}$. $\hfill\Box$

\vspace{1cm}

\noindent
{\it Remark:} Suppose, as it is the case for $\bpsi$, that  we only know that $\ds\int \rho |\bpsi_1-\bpsi_2|^p \leq \epsilon ^p$, then we have instead of (\ref{2M}),
\be
\epsilon^p\geq\int_{B_r}\rho|\bpsi_1-\bpsi_2|^p \geq \rho(B_{M/2R}(x))  M^ {p+d}.
\en 

\vspace{1cm}
\noindent
The first part of the lemma yields immediately that $\bfi \in C^{\alpha}(I, C^0(\Omega))$ for some $\alpha>0$.
Moreover if $\phi^*_1, \phi^*_2$ are the Legendre transform of $\phi_1, \phi_2$, then $\|\phi_1^*-\phi_2^*\|_{\Linf}
 \leq \|\phi_1-\phi_2\|_{\Linf}$, thus $\bpsi \in C^{\alpha}(I, C^0(B_R))$, and
this gives the point 2.
\\
The second point of the lemma will be used to prove point 4 and 5:
Indeed, if $\rho$ supported in $\bar\Omega'$ for some open set $\Omega'$, and there exists $0<\lambda, \Lambda$ such that $\lambda \leq \rho \leq \Lambda$ in $\Omega'$,
from Theorem \ref{2caf} we get that for any $\omega' \subset\subset \Omega'$,
$\bpsi(t,\cdot) \in C^{1,\alpha_1}(\omega')$ for some $\alpha_1>0$.
Since $\dt \bfi \in L^{1*}(\Omega)$, using (\ref{2identite2}) we 
get that 
\be
\int \rho_n|\partial_t\bpsi_n|^{1*} \leq C
\en
uniformly in $n$, and thus that 
\be
\partial_t\bpsi_n \in \Linf(I, L^{1*}(\Omega')).
\en
Therefore we can use Lemma \ref{2Calpha} to obtain that for any $\omega' \subset\subset \Omega'$, 
$\ds \nabla\bpsi\in C^{\beta}(I, C^0(w'))$ 
(point 4 of Theorem \ref{2main1}).
\\
Under the additional assumption that $\Omega'$ is convex, Theorem \ref{2caf}
yields that $\bfi(t,\cdot)$ in $C^{1,\alpha_2}(\Omega)$ for some $\alpha_2>0$. The same procedure
as above yields point 5.

\vspace{1cm}
\noindent
Now we prove the point 3 of Theorem \ref{2main1}:

\subsubsection*{Proof of the bound on $\dt g$}

Recall from Theorem \ref{2main2}:
\be
\int_{\Rd} \rho_n |\partial_t\nabla\bpsi_n| \leq C(d,R)\|\rho_n\|^{\demi}_{\Linf(B_R)}\|\rho_n v_n^2\|^{\demi}_{\Linf(B_R)}
\en
We have 
$g(t,a)=\nabla\bpsi(t,{\bf X}(t,a))$ and thus formally 
$$\partial_t g(t,a)= \partial_t\nabla\bpsi(t,{\bf X}(t,a))+ D^2\bpsi(t,{\bf X}(t,a))\partial_t {\bf X}(t,a).$$
Since $\rho_n$ converges strongly (actually weakly would be enough) to $\rho$,
we know that $\nabla\bpsi_n$ converges almost everywhere to $\nabla\bpsi$. (See \cite{Br1} for a 
proof of this fact,  which relies on the convexity of $\bpsi_n$ and on the uniqueness of the polar factorization). 
Now consider  
\be
g_n(t,a)=\int_{\Rd} \nabla\bpsi_n(t,y)\eta_n(y-{\bf X}(t,a))dy=(\eta_n*\nabla\bpsi_n) (t,{\bf X}(t,a))
\en 
with $\eta_n$ a smoothing kernel as above.
Then $g_n$ converges almost everywhere to $g$.
For $f \in C^0(I\times\Omega, R^d)$, let us compute
\be
\int_I\int_{\Omega}\partial_t  g_n(t,a)\cdot f(t,a)\,dtda=T_1+T_2,
\en
with
\be
T_1 &=&\int_I\int_\Omega\int_{\Rd} \eta_n(y-{\bf X}(t,a))\partial_t\nabla\bpsi_n(t,y)\cdot f(t,a)\, dy da dt\\
T_2 &=& - \int_I\int_\Omega\int_{\Rd}\nabla\bpsi_n(t,y)\cdot f(t,a) \  \partial_t {\bf X}(T,a)\cdot \nabla\eta_n(y-{\bf X}(t,a))\, dy da dt
\en
Let us evaluate $T_1$ and $T_2$.
\be
|T_1|&\leq &\int_I \|f(t,.)\|_{\Linf(\Omega)}\int_{\Rd\times\Rd}\rho(x)\eta_n(y-x)|\partial_t\nabla\bpsi_n(t,y)|\,dxdy dt\\
&\leq& \int_I \|f(t,.)\|_{\Linf(\Omega)}d_n \int_{\Rd\times\Rd}\rho_n(y)|\partial_t\nabla\bpsi_n(t,y)|\,dxdy\\
&\leq& \int_I \|f(t,.)\|_{\Linf(\Omega)}C(R,d)\|\rho_n\|_{\Linf(I\times \Rd)}\|v_n\|_{\Linf(I\times \Rd)}
\en
with $d_n=1/c_n$ and from Theorem \ref{2main2}.
For $T_2$ we have:
\be
|T_2|&=&\left|\int_I\int_\Omega\int_{\Rd}\nabla\bpsi_n(t,y)\cdot f(t,a) \ \partial_t {\bf X}(T,a)\cdot \nabla\eta_n(y-{\bf X}(t,a))\, dy da dt\right|\\
&=&\left|\int_I\int_\Omega\int_{\Rd} \partial_t {\bf X}^t(T,a)\cdot(D^2\bpsi_{n}*\eta_n)(t,{\bf X}(t,a))\cdot f(t,a)\, dy da dt\right|\\
&\leq &  \int_I \|f(t,.) \partial_t {\bf X}(t,.)\|_{\Linf(\Omega)} \int_{\Rd} \rho (t,x)(|D^2\bpsi_{n}|*\eta_n)(x)dx \ dt\\
&\leq & \int_I \|f(t,.) \partial_t {\bf X}(t,.)\|_{\Linf(\Omega)} \|\rho(t,.)\|_{\Linf(\Rd)}C(R,d,\Omega) \ dt
\en
where we have used the bound on $\|D^2\bpsi\|_{L^1_{loc}}$ (Lemma \ref{2D2phiL1});
we conclude that
\be
\|\partial_t g\|_{\Linf(I, {\cal M}(\Omega))}\leq C(R,d,\Omega)\|\rho\|_{\Linf(I\times B_R)}
\|\partial_t {\bf X}\|_{\Linf(I\times \Omega)}.
\en
This achieves the proof of Theorem \ref{2main1}.

$\hfill\Box$


\section{Proof of Theorem \ref{2main2}}

\subsection{H\"older regularity}
It has been established ((\ref{2diverg}) and Theorem \ref{2regulier}) that $\dt\bpsi_n$ satisfies
\be
\nabla\cdot(M_n \dt\nabla\bpsi_n)=\sum_{i,j}M_{n,ij}\partial_{ij}\partial_t\bpsi_n= -\dt\rho_n=  -\nabla\cdot(\rho_n v_n)
\en
where $M_n$ is the comatrix of $D^2\bpsi_n$.
To establish the H\"older regularity of $\dt\bpsi_n$ we need to combine three preliminary results: 

The first one (Theorem \ref{2theocagu})
asserts the Harnack inequality for solutions of the homogeneous linearized Monge-Amp\`ere equation under a condition which is satisfied when the density $\rho$ is between two positive constants. 

The second 
one (Theorem \ref{2main5}, Theorem \ref{2main4} and Corollary \ref{2corlocal}) is a local maximum principle that
generalizes the local maximum principle for uniformly elliptic equations, to degenerate elliptic equations of the 
form $\nabla\cdot (M\nabla u)= \nabla\cdot f$.  The uniform ellipticity is relaxed to the condition that the (positive symmetric matrix valued) 
functions $M, M^{-1}$ belong to $L^p$ for $p$ large enough. $p$  depends only on the dimension $d$.

The third one (Theorem \ref{2W2p}) asserts that the comatrix of $D^2\bpsi$, and its inverse, are indeed in $L^p_{loc}(\Omega')$
 provided that the density $\rho$ is close enough to a continuous positive function, the closeness being measured in $\Linf$ norm.

\bigskip

The result will be a consequence of the following propositions:

\begin{prop}\label{2main3c}
Let $\rho=\bx\#da$ be supported in $\bar\Omega'$, $\lambda$ and $\Lambda$ be two positive constants such that $0<\lambda\leq \rho (t,x) \leq \Lambda$ for all $(t,x) \in I\times\Omega'$. Let $\rho_n, v_n$ be constructed from $\bx$ as above. Let $(\bfi_n, \bpsi_n)$ be associated to $(\rho_n, \Omega)$ through (\ref{2defPhi}, \ref{2defPsi}).
Then there exists $\beta<1$, and  for any $p>d$, there exists $C$ such that
for any $\omega' \subset\subset\Omega'$, for any $(y,r)$ with $ B_{4r}(y)\subset \omega'$, 
\be
osc_{\dt\bpsi_n}(r/2) \leq \beta osc_{\dt\bpsi_n}(r) + C r^{\delta}
\en
for $n$ large enough.
$\beta<1$ depends on $(\lambda, \Lambda)$ (see Theorem \ref{2theocagu}),  
$C$ depend on $(p,\lambda, \Lambda,\inf_{x\in \omega'} d(x,\partial\Omega')$, $\|D^2\bpsi\|_{L^p(B_r(y))})$,
 $\ds\delta=1-\frac{d}{p}$, and  
\be
osc_u(r)=\max_{B_r}u-\min_{B_r}u.
\en
\end{prop} 
{\it Remark:} The requirement $n$ large enough is just to enforce that $\lambda\leq \rho_n \leq \Lambda$.

\vspace{1cm}

\begin{prop}\label{2main3d}
Under the assumptions of Proposition \ref{2main3c}, we have, for every $\omega'\subset\subset \Omega'$, 
\be
\|\dt\bpsi_n\|_{\Linf(\omega')} \leq C(K, p, \inf_{x\in \omega'} d(x,\partial\Omega'), \lambda, \Lambda, \|v_n\|_{\Linf(d\rho(t))})
\en 
where $K=\|D^2\bpsi_n+ D^2\bpsi^{-1}_n\|_{L^p(\omega')}$, $p>d$.
\end{prop}

\vspace{1cm}

\begin{prop}\label{2renorm}
Under the assumptions of Proposition \ref{2main3c}, for any $p < \infty$, there exists $\epsilon>0$
such that if $|\rho-1|\leq \epsilon$ in $\Omega'$, then for every $K'\subset \Omega'$, $K'$ compact, there exists $C_{K'}$ such that
\be
\limsup_n\|D^2\bpsi_n + D^2\bpsi_n^{-1}\|_{L^p(K')} \leq C_{K'}.
\en
\end{prop}
\vspace{1cm}
\noindent
Temporarily admitting these propositions we obtain the following:

\subsubsection*{Proof of the first part of Theorem \ref{2main2}}
From Propositions \ref{2main3c}, \ref{2main3d}, \ref{2renorm},
we obtain that for any $\omega'\subset\subset \Omega'$, there exists $C_{\omega'}$, $\beta<1$ independent of $n$  such that,
for $n$ large enough, for any $B_r=B_r(y)\subset\omega'$, with $B_{4r}\subset \Omega'$, we have:
\be
osc_{\dt\bpsi_n}(r/2) \leq \beta osc_{\dt\bpsi_n}(r) + C_{\omega'} r^{\delta}.
\en
Moreover from Proposition \ref{2main3d}, $\dt\Psi_n$ is uniformly bounded for the sup norm inside $\omega'$.
It is well known that this property implies H\"older continuity: using \cite{GT}, Lemma 8.23, 
we obtain that for $n$ large enough, 
for any $\omega'\subset\subset \Omega'$,  there exists $\alpha >0, C_{\omega'}$ that do not depend on $n$, such that for any $(x,y)\in \omega'$,
\be
|\dt\bpsi_n(y)-\dt\bpsi_n(x)|\leq C_{\omega'} |x-y|^{\alpha}.
\en
Thus we have a uniform $\Linf(I;C^{\alpha}(\omega'))$ bound 
that will pass to the limit as $n\to \infty$. We thus obtain the $C^{\alpha}$ estimate of Theorem \ref{2main2}.
\\
To obtain H\"older continuity for $\dt\bfi$, in the case where $\Omega'$ is convex,
we just have to use the identity (\ref{2identite0}) 
\be
\dt \bfi = - \dt \bpsi(\nabla\bfi)
\en
and the H\"older regularity of $\nabla\bfi$, under the condition 
$0< \lambda\leq \rho\leq\Lambda$, $\Omega'$ convex (Theorem \ref{2caf}), 
to conclude H\"older regularity for $\dt \bfi$.

$\hfill \Box$

\vspace{1cm}
\noindent
In the next proofs we drop the suffix $n$ for simplicity.
\\
{\it Proof of Proposition \ref{2main3d}:}
This proposition is a direct consequence of Theorem \ref{2main4}.
It has been established that $\dt\bpsi$ satisfies
\be
\nabla\cdot(M \dt\nabla\bpsi)=\sum_{i,j}M_{ij}\partial_{ij}\partial_t\bpsi= -\dt\rho=  -\nabla\cdot(\rho v)
\en
where $M$ is the comatrix of $D^2\bpsi$, given by $M=\det D^2\bpsi [D^2\bpsi]^{-1}$ or $M=\rho D^2\bfi(\nabla\bpsi)$.
We remember that $0<\lambda\leq \rho\leq \Lambda$. From Theorem \ref{2main1}, we have the a priori bound  
\be
\int_{\Omega} |\dt\bfi|^{1*}\leq C (\|\rho_n |v_n|^2\|_{\Linf}, \Omega, R, d).
\en
Using then that $\dt \bpsi= -\dt\bpsi(\nabla\bfi)$ we have
\be
\int \rho|\dt\bpsi|^{1*} = \int_{\Omega}|\dt \bfi|^{1*}
\en
and thus 
\be
\int_{\Omega'} |\dt\bpsi|^{1*}\leq \frac{C}{\lambda}.
\en
We can therefore apply Theorem \ref{2main4} with $a_0=1*$.

$\hfill\Box$

\vspace{1cm}
\noindent
{\it Proof of Proposition \ref{2main3c}:}
\\
We consider a ball $B_{4r}(y)$ contained in $\Omega$ and write $\dt \bpsi= u + w$ where
$u$ satisfies
\be
&&\nabla\cdot (M \nabla u) = -\nabla\cdot ( \rho v),\\
&&u=0 \text{ on } \partial B_r(y),
\en
\\
and $w$ satisfies
\be
&&\nabla\cdot (M \nabla w)=0\\
&& w=\dt \bpsi \text{ on } \partial B_r(y).
\en
Note that $w$ satisfies also $M_{ij}\partial_{ij}w = 0$ which is the equation treated in
\cite{CaGu}.
\\
We denote $\ds osc_f(r)=\sup_{B_r}f - \inf_{B_r} f$ and 
$\ds osc_f(\partial B_r)= \sup_{\partial B_r}f - \inf_{\partial B_r} f$.
\\
The assumptions of Theorem \ref{2theocagu} are satisfied: indeed,
in $\omega' \subset \subset \Omega'$, we have, for $n$ large enough,
$\lambda \leq \rho_n \leq \Lambda$.
From Theorem \ref{2theocagu}, there exists $\beta<1$ such that
\be
&& osc_w(r/2)\leq \beta osc_w(r).
\en 
From Corollary \ref{2corlocal} we have 
\be
&& \sup_{B_{r}} |u|  \leq C \|\rho v\|_{\Linf}r^{\alpha},
\en
where $\alpha=1-d/p$, $\ds C= C_0\|M^{-1}\|_{L^p}=C_0\|\rho^{-1}D^2\bpsi\|_{L^p(B_{r})}$ (note that we have $0<\lambda\leq \rho \leq \Lambda$). 
Combining the two estimates, we have
\be
osc_{\dt\bpsi}(r/2) &\leq& osc_{w}(r/2) + osc_{u}(r/2)\\
&\leq& \beta osc_w(r)+Cr^{\alpha}\\
&\leq& \beta osc_w(\partial B_r) + Cr^{\alpha}\\
&\leq& \beta osc_{\dt\bpsi}(\partial B_r) + Cr^{\alpha}\\
&\leq& \beta osc_{\dt\bpsi}(r) + Cr^{\alpha}
\en
where in the third line we have used the maximum principle to say that
$osc_w(r)=osc_w(\partial B_r)$ since $w$ can not have interior extrema.
Finally we conclude
\be
&& osc_{\dt \bpsi}(r/2)\leq \beta osc_{\dt \bpsi}(r) + C r^{\alpha}.
\en
This achieves the proof of Proposition \ref{2main3c}.
$\hfill\Box$

\vspace{1cm}
\noindent
{\it Proof of Proposition \ref{2renorm}}
We show here how to use the $W^{2,p}$ regularity Theorem \ref{2W2p} to obtain estimates. 
First let us notice that 
if $\nabla \bpsi$ satisfies (\ref{2defPsi}) for  $\rho$ supported in $\bar\Omega'$, $0<\lambda \leq\rho\leq \Lambda$, and  since $\Omega$ is convex,  we know from \cite{Ca3} that
$\bpsi$ is strictly convex in $\Omega'$ and solution in the viscosity sense to 
\be
\det D^2\bpsi=\rho
\en 
in $\Omega'$. Moreover $\bpsi$ is $C^{1,\alpha}_{loc}$ in $\Omega$ (Theorem \ref{2caf}).
From the strict convexity, for any $x\in \Omega'$, there exists 
a section 
\be
S_{t_x,x}=\{y: \bpsi(y)\leq \bpsi(x)+\nabla\bpsi(x).(y-x)+t_x\} 
\en
with non-empty interior and compactly contained in $\Omega'$. (Indeed the strict convexity means that diameter of the sections decreases  to 0 as
the height of the section $t_x$ goes to 0).
Then for every compact set $K$ contained in $\Omega'$ there
exists a finite covering of $K$ by sets $\frac{1}{3d}S_i$, $S_i=S_{t_{x_i},x_i}$, and $\frac{1}{3d}S_i$ means a contraction of $S_i$ with respect to $x_i$.
Then the functions $u_i(y)=\bpsi(y)-t_i - \nabla\bpsi(x_i)\cdot(y-x_i)$
are solutions of 
\be
&&\det D^2 u_i = \rho \text{ in } S_i\\
&&u_i=0\text{ on } \partial S_i.
\en
From John's lemma (see \cite{Ca0}), we can find an affine transformation $T_i$, with $\det T_i=1$
and a real number $\mu_i$
such that $B_{1} \subset \mu_i^{-1} T_i^{-1}(S_i)=\tilde S_i \subset d B_{1}$. Finally,
considering $\tilde u_i(y) =\frac{1}{\mu_i^2}u_i(\mu_i T_i \  y)$ 
we get that $\tilde u_i$ is solution to
\be
&&\det D^2 \tilde u_i(y) = \tilde\rho(y)=\rho(\mu_i T_i y) \text{ in } \tilde S_i\\
&&\tilde u_i=0\text{ on } \partial \tilde S_i\\
&&B_1(x_i)\subset \tilde S_i \subset d B_1(x_i).
\en 
We can invoke Theorem \ref{2W2p} for $\tilde u_i$:
For any $0<p<\infty$, if $|\tilde\rho-1|\leq \epsilon(p)$ (this property is invariant under the renormalizations
performed above), we have
\be
&&\|D^2\tilde u_i + D^2\tilde u_i ^{-1}  \|_{L^{p}(B_{\demi})}\leq C\\
&&\left( \text{meas}(S_i) \right)^{-1/p}\|D^2 u_i+ D^2 u_i^{-1}\|_{L^{p}(\frac{1}{2d}S_i)}\leq C\|T_i\|^2.
\en
By our covering process, we have $K\subset \bigcup_i T_i \mu_i B_{\frac{1}{3}}(x_i)$. It follows that for every compact set $K\subset \Omega'$, there exists and constant $C_{K}$ such that $\|D^2\bpsi\|_{L^p(K)}\leq C_K$ 
and $\|D^2\bpsi^{-1}\|_{L^p(K)}\leq C_K$. The constant $C_K$ depends on the  supremum of the norm of the transformations $T_i$
and can be taken (by compactness) uniformly bounded given $\Omega, \Omega', K, \lambda, \Lambda$.

Now we show that this covering process behaves uniformly well when we consider the regularization $\rho_n$ of $\rho$ and let               
 $n$ go to $\infty$. Indeed the corresponding $\bpsi_n$ will converge uniformly to $\bpsi$ and since the limit $\bpsi_n$
is $C^1$ the sequence $\nabla\bpsi_n$ converges also uniformly in every compact set of $\Omega'$.
Therefore the set $S_{i}^n=\{y, \bpsi_n(y) \leq \bpsi_n(x_i) + \nabla\bpsi_n(x_i)\cdot(y-x_i) + t_i\}$ converge uniformly to $S_i$.  This means that for $n$ large enough, the set $K$ will be covered by $\bigcup_i \frac{1}{2d}S_i^n$.
Consider $\mu_i^n, T_i^n$ the corresponding normalization.
then we also have $T_i^n, \mu_i^n$ converging  to $T_i, \mu_i$, and $K$ will be covered by $\bigcup_iT_i^n \mu_i^n B_{\frac{1}{2}}(x_i)$.

Moreover since we consider a compact set $K$ contained in $\Omega'$ and since $|\rho-1| \leq \epsilon$ in $\Omega'$, it follows from the construction of $\rho_n$ that, for $n$ large enough, $|\rho_n-1|\leq \epsilon$ in $\Omega'$. 
For $n$ large enough, the functions $\tilde u_i^n$ (obtained by the renormalization procedure) will thus all satisfy the assumptions of Theorem \ref{2W2p}.
\\
Therefore, for every $K\subset\subset \Omega'$,  there exists $C_K$ independent of $n$ such that, for $n$ large enough,
\be
\|D^2\bpsi_n + D^2\bpsi_n^{-1}\|_{L^p(K)} \leq C_K.
\en
This achieves the proof of Proposition \ref{2renorm}.

$\hfill \Box$

\subsubsection*{Proof of the gradient bounds}
This is point 2  of Theorem \ref{2main2}.
The gradient bounds follow directly from Proposition \ref{2main3b} combined with Proposition \ref{2renorm}. In estimates (\ref{2dtgradphi}, \ref{2dtgradpsi}) take
$r=\infty$. Note that from Lemma \ref{2D2phiL1} we have
the bound $\|D^2\bpsi\|_{L^1(B_R)}\leq C(R,d,\Omega)$.
This ends the proof of Theorem \ref{2main2}.

$\hfill \Box$

\section{The periodic case: proof of Theorem \ref{2main2per}}
This result is only an adaptation of the two previous Theorems.
All the regularity results used adapt to the periodic case as follows:
\begin{theo}\label{2theoperio}
Let  $\rho$ be a Lebesgue integrable probability measure on $\Rd/\Zd$. There exists a
unique $\bpsi$ convex on $\Rd$, with  $\bpsi - |x|^2/2$ periodic, that satisfies \be
\det D^2\bpsi = \rho
\en
 in the following sense:
\be
\forall f \in C^0(\Rd/\Zd),  \ \int_{\Td} \rho f(\nabla\bpsi)= \int_{\Td} f.
\en
It has the following regularity properties:
\begin{enumerate}
\item If for some pair $(\lambda, \Lambda)\in \R_+^*$, we have $\lambda\leq \rho \leq \Lambda$,  then for some $\alpha>0$ depending on $\Lambda/\lambda$, $\bpsi-|x|^2/2$ is in $C^{1,\alpha}(\Td)$.
\item For every $p<\infty$, there exists $\epsilon(p)$ such that if $|\rho-1|\leq \epsilon(p)$, then $\bpsi-|x|^2/2 \in W^{2,p}(\Td)$.
\item If $\rho$ is positive and in $C^{\infty}(\Td)$, then $\bpsi -|x|^2/2 \in  C^{\infty}(\Td)$.
\end{enumerate}
\end{theo} 
We then modify the approximation procedure as follows: we take 
\be
&&\rho_n = c_n(\eta_n* \rho + \frac{1}{n})\\
&&\det D^2\bpsi_n = \rho_n
\en  
with the constant $c_n$ such that $\int_{\Td} \rho_n = 1$.
Then we use the same techniques as in the Theorems \ref{2main1}, \ref{2main2}. 
\\
We only mention the two new results that arise in this case:

In point 2, we obtain that $\bg \in C^{\alpha}(I, \Linf(\Td))$. Indeed, 
$\bg= \nabla\bpsi(t,\bx(t))$. We already know that, under the present assumptions, $\nabla\bpsi \in C^{\alpha}(I \times \Td)$, moreover $\bx \in W^{1,\infty}(I, \Linf(\Td))$  and the result follows.

In point 4, under the assumption that $\|\rho-1\|_{\Linf(I\times\Td)} \leq \epsilon$ for $\epsilon$ small enough depending on $q$, ,
we are able to obtain a bound in $L^q(\Td),  \ q<2$ for $\dt \bg$. Indeed, writing 
\be
\bg_n(t,a)=\nabla\bpsi_n(t,X(t,a))
\en 
as in the proof of Theorem \ref{2main1}, and differentiating with respect to time, we obtain
\be
\dt\bg_n(t,a)=\dt\nabla\bpsi_n(t,\bx(t,a))+ D^2\bpsi_n(t,\bx(t,a)) \dt\bx(t,a).
\en
with $\bpsi_n$ obtained from $\rho_n$, and thus in $C^{\infty}(I\times\Td)$.
If $\rho$ is close enough to 1 so that $D^2\bpsi_n$ is bounded in $L^p(\Td)$ (cf. Theorem \ref{2theoperio} above), 
the first term is bounded in $L^q(\Td)$, with $q=\frac{2p}{1+p}$ (as in Proposition \ref{2main3a}). The second term is bounded in $L^p(\Td)$. Then we let $\bg_n$ converge to $\bg$.

Note that this bound can not be obtained in the non periodic case since we have only interior regularity available for $\bpsi$.

\section{Counter-examples}\label{contrex}
Here we show through some examples that the bounds  obtained in Theorem \ref{2main1} are sharp under 
our present assumptions.

{\bf Example 1:} $\dt\nabla\bfi \notin L^1_{loc}$ and $\dt\bfi \notin C^0$. 
\\
Consider in $\Omega=B(0,1)$ in $\R^2$,
and ${\bf X}(t,\cdot): B(0,1)\to \R^2$ defined with complex notations ${\bf X}=x + i y$ by

on $y>0$, 
\be
{\bf X}(t,(x,y))=e^{it} (x+iy) + it,
\en

on $y<0$, 
\be
{\bf X}(t,(x,y))=e^{it} (x+iy) + t^2.
\en
We check that ${\bf X}\# da$ has a density bounded by 1, that $\dt {\bf X} \in \Linf(\Omega\times \R^+)$. 
If ${\bf X}=\nabla\bfi\circ g$ is the polar factorization of ${\bf X}$
then up to a constant, $\bfi$ is defined for $t>0, (x,y)\in \Omega$ by:
\be
\bfi(t,(x,y))= \sup\{\demi(x^2 +y^2) + t^2 x, \demi(x^2 +y^2) + t y \}.
\en
On $\{ y > tx\}$  we have
\be
&&\bfi(t,(x,y)) = \demi(x^2 +y^2) + t y,\\
&& \nabla\bfi(t,(x,y))=(x,y) + (0,t), 
\en
and on $\{ y < tx\}$ 
\be
&&\bfi(t,(x,y)) = \demi(x^2 +y^2) + t^2 x,\\
&& \nabla\bfi(t,(x,y))=(x,y) + (t^2,0).
\en
Thus 
\be
&&\dt \bfi (t,(x,y))= y \chi_{\{ y > tx\}}+ 2tx  \chi_{\{ y < tx\}}  \notin C^0,\\
&&\dt\nabla\bfi(t,(x,y))= (0,1)\chi_{\{ y > tx\}}+ (2t,0)\chi_{\{ y , tx\}} +
(t^2,-t) {\cal H}^{d-1}\{y=tx\} \notin L^1_{loc}
\en

{\bf Example 2:} Here we adapt a counterexample of Wang to build
an example of a solution where $\dt \bpsi \notin C^0$.
\\
In $\Rd$, let $x=(x_i)_{1\leq i\leq d}$ and 
\be
{\bf X}(0,x)= \nabla\bfi_0(x)
\en
$\bfi_0(x)$ convex Lipschitz on $\Omega$,  $\bfi=+\infty$ outside, such that
$\rho=\nabla\bfi_0(x)\# dx$ has a density in $\Linf(\R^2)$. 
Let 
\be
{\bf X}(t,x)= \nabla\bfi_0(x) + t v
\en
for some fixed $v\in \Rd$. ${\bf X}$ is Lipschitz with respect to time.
Then 
\be
&&\bfi(t,x)=\bfi(x)+ t x\cdot v,\\
&&\nabla\bfi(t,x)= \nabla\bfi_0(x) + t v.
\en
If $\bpsi_0$ is the Legendre transform of $\bfi_0$, the Legendre transform of $\bfi(t,\cdot)$ is given by  
\be
&&\bpsi(t,x)=\bpsi_0(x-tv),\\
&&\nabla \bpsi(t,x)=\nabla \bpsi_0(x-tv),
\en 
thus 
\be
&&\partial_t\bpsi(t,x)=v\cdot \nabla\bpsi_0(x-tv),\\
&&\dt\nabla\bpsi(t,x)=D^2\bpsi_0(x-tv)\cdot  v.
\en
Wang has shown in \cite{W} some counterexamples to the regularity of solutions of Monge-Amp\`ere equations: namely, for $d\geq 3$ he has exhibited a solution $u$
of
\be
\det D^2 u= f
\en with $f$ only bounded by above, such that $u\notin C^1$.
By taking $\bpsi_0= u$ one has an example of time dependent map such that 
\be
\partial_t\bpsi(t,x)=v\cdot \nabla\bpsi_0(x-tv) \notin C^0.
\en

\newpage

\section{Application: the semi-geostrophic
equations}\label{2SG} The semi-geostrophic system is derived as an approximation to the primitive equations in meteorology, and is believed  to model frontogenesis (see \cite{CP}). 
The formulation of the 3-d incompressible version is the following:
we look for a time dependent probability measure $\rho$ that satisfies the following $SG$ system:
\beq 
&&\dt \rho +\nabla\cdot (\rho \bv)=0\label{2sg1}\\ 
&&\bv(t,x)=\left(\nabla\bpsi(t,x)-x\right)^{\perp}\label{2sg2}\\ 
&&\det D^2\bpsi(t,x)= \rho(t,x)\label{2sg3}. 
\enq
Here ${\bf v}^{\perp}$ means $(-{\bf v}_2, {\bf v}_1, 0)$.
Equation (\ref{2sg3}) is understood in the sense of (\ref{2defPsi}), where an open set $\Omega$ of total mass 1 has been given before. 
\\
The system has also a periodic version in which $\Omega=\T^3$ itself and equation (\ref{2sg3}) is solved with the condition that $\bpsi-|x|^2/2$ is $\Z^3$ periodic.
\\
The set $\Omega$ is here called the physical space, whereas the space in which $\rho$ lives is the dual space.
Existence of global weak solutions for the $SG$ system with initial data in $L^1$ has been proved in \cite{BB2}, \cite{CuG} and \cite{LN}.
Note that uniqueness of weak solutions is still an open question.

\subsection{The Lagrangian formulation of the $(SG)$  system}
Here we look for a mapping $\bx: \R^+ \times \Omega \to \R^3$ that satisfies
\beq
&&\dt {\bf X}(t,a) = \left(\nabla\bpsi(t,{\bf X}(t,a)) -{\bf X}(t,a)\right)^{\perp}\label{2trajec}\\
&&\nabla\bpsi(t) \circ \bx(t)=\bg(t) \in G(\Omega), \ \bpsi \text{ convex}\label{2trajec1}. 
\enq 
If we define $\rho(t)=\bx(t)\#da$, the last equation means that for all $t$, $\bpsi(t)$ solves
$\det D^2\bpsi(t)=\rho(t)$ in the sense of (\ref{2defPsi}). Having $\bx$ solution of (\ref{2trajec}, \ref{2trajec1}) implies that $\rho(t)=\bx(t)\#da$ is solution of (\ref{2sg1}, \ref{2sg2}, \ref{2sg3}). $\bx$ defines the characteristics in the dual space whereas $\bg$ defines the characteristics in the physical space.

We expose briefly the arguments that allow to define the characteristics of the $SG$ system:
\\
1- First we check that $\bx(t)$ will satisfy for any time $t$ the condition (\ref{2N}): indeed, the flow being incompressible, all the $L^p$ norms of $\rho$ are conserved. Therefore, given the potential $\bpsi(t)$, if $X_0$ satisfies the condition (\ref{2N}), or equivalently if $\rho_0\in L^1$, 
then we know a priori that $\bx(t)$ satisfies the condition (\ref{2N}) for all time. 
\\ 
2- The velocity field is a priori bounded in $BV$ because of the convexity of $\bpsi$ (see Lemma \ref{2D2phiL1}). Moreover it is incompressible.
Therefore thanks to the  result of \cite{Am}, the characteristics of the corresponding ODE are 
uniquely defined for almost every initial data,
which means that the curve $t\mapsto \bx(t,a)$ is uniquely defined for almost every $a\in \Omega$.

For $\Omega$ bounded, it is easily checked (see  \cite{BB2}) that if $\bx_0 \in \Linf(\Omega)$, then  $(\bx,\dt\bx) \in \Linf([0,T]\times\Omega)$ for all $T>0$.
The velocity field being incompressible, if $\rho_0\in \Linf(\R^3)$, then 
$\rho\in \Linf(\R^+\times\R^3)$.
Note that the Lagrangian system can also be defined in a periodic space, where $\bx$ is periodic in space for all time, and we require $\bpsi -|x|^2/2$ to be periodic.
The bound of $\bx, \dt\bx$ in $\Linf(\R^+\times \T^3)$ is then independent of the initial data.
Moreover, in this setting, if  $\rho_0$ is such that 
\beq\label{2semigeovide}
0<\lambda \leq \rho_0 \leq \Lambda 
\enq 
for two constants $\lambda, \Lambda$,
this property remains satisfied for all time, once again due to the incompressibility of the velocity field.
\\
Thus we conclude the following:
\begin{lemme}
Let $\bx_0\in \Linf(\Omega;\R^3)$, $\rho_0=\bx_0\#da\in \Linf(\R^3)$.
Then $\rho, \bx$ the corresponding solution of the $SG$ system
satisfies for all $T>0$, 
\be
&&\bx, \dt\bx \in \Linf([0,T]\times \Omega)\\
&&\rho \in \Linf(\R^+\times \R^3).
\en
In the periodic case this remains true,  and if moreover $\rho_0$ satisfies (\ref{2semigeovide}), then for all time $t$, $\rho(t)$ satisfies (\ref{2semigeovide}). 
\end{lemme}
Under the assumptions of the above lemma, it is clear that $\bx$ satisfies the assumptions of Theorem \ref{2main1}.  
In the periodic case, if satisfied at time 0, all the assumptions of Theorem \ref{2main2per} are satisfied for all time. 
We can now state the following theorem of partial regularity. We restrict ourselves to the periodic case.
\\
{\it Remark:} We also conjecture that the assumptions of Theorem \ref{2main2} can be satisfied for some finite time, 
but the control the evolution of the support of $\rho$ poses some some difficulties.

\begin{theo}
Let $\bx, \rho, \bg, \bpsi,\bfi$ be as above, with $\rho=\bx\#da$ be a space-periodic solution of (\ref{2sg1}, \ref{2sg2}, \ref{2sg3}), and $\bx$ the corresponding space-periodic solution of (\ref{2trajec}, \ref{2trajec1}).  
Suppose that $\rho_0 \in \Linf(\T^3)$, then 
\be
&&\dt \bg \in \Linf(\R^+, {\cal M}(\T^3)), \\ 
&&\dt \nabla\bfi \in \Linf(\R^+, {\cal M}(\T^3)).
\en
If moreover there exists $0<\lambda,\Lambda$ such that $\lambda\leq \rho_0\leq \Lambda$, then there exists $\alpha>0$ depending on $(\lambda, \Lambda)$ such that
\be
\bg \in C^{\alpha}(\R^+, L^{\infty}(\T^3)).
\en
For all $p<2$, there exists $\epsilon(p)$, such that
if $|\rho_0-1|\leq  \epsilon(p)$, then 
\be
\dt \bg \in \Linf([0,T], L^p(\T^3)).
\en
There exists $\epsilon_0$, such that if $|\rho_0-1| \leq \epsilon<\epsilon_0$, then
\be
\dt \bfi, \dt\bpsi \in \Linf(\R^+, C^{\alpha}(\T^3))
\en
where $\alpha>0$ depends on $\epsilon$.
\end{theo}

{\bf Remark: The equations of motion in physical space}
We derive  here formally the equation giving the evolution of $\bg$: writing (\ref{2hodge-ampere}) with $\bv$ as above, we have
\be
&&(x-\nabla\bfi)^{\perp}=\bv(\nabla\bfi)=\dt \nabla\bfi + D^2\bfi w,\\
&&\nabla\cdot w =0,
\en
where $\dt\bg(\bg^{-1})=w$. 
This equation formally determines the evolution of the system, since the knowledge of $\bfi(t)$ determines a unique pair $\dt\nabla\bfi,w$ satisfying the above decomposition (see Proposition \ref{2hodge}). One can see a parallel with the Euler incompressible equation
where the evolution is given by solving the following decomposition problem:
\be
&& - v\cdot\nabla v = \dt v + \nabla p,\\
&& \nabla\cdot v = 0.
\en 
Thus 
the semi-geostrophic equations are associated to the decomposition of vector fields of Proposition \ref{2hodge}
in a similar way
as the Euler incompressible equations are associated to the Hodge ``div-curl'' decomposition. 

\newpage

\bibliography{biblio}

\end{document}